\documentclass[11pt,a4paper]{article}
\usepackage{ifthen,latexsym,amssymb,amsmath,bbm,fixmath}
\usepackage[nobysame,initials]{amsrefs}
\usepackage{dsfont}
\usepackage{authblk}

\newcommand{\C}[1]{{\protect\mathcal{#1}}}

\newcommand{\I}[1]{{\mathbb #1}}
\renewcommand{\O}[1]{\overline{#1}}

\newcommand{\V}[1]{\mathbold{#1}}
\newcommand{\nib}[1]{\noindent{\textbf{#1}}}

\newcommand{\e}{\varepsilon}

\renewcommand{\mid}{:}
\renewcommand{\ldots}{\hspace{0.9pt}.\hspace{0.3pt}.\hspace{0.3pt}.\hspace{1.5pt}}
\renewcommand{\ge}{\geqslant}
\renewcommand{\le}{\leqslant}

\renewcommand{\succeq}{\succcurlyeq}

\newif\ifnotesw\noteswtrue

\noteswfalse	
\newcommand{\hide}[1]{}

\newcommand{\beq}[1]{\begin{equation}\label{#1}}
\newcommand{\eeq}{\end{equation}}

\newtheorem{theorem}{Theorem}[section]
\newtheorem{lemma}[theorem]{Lemma}

\newtheorem{proposition}[theorem]{Proposition}

\newtheorem{conjecture}[theorem]{Conjecture}
\newtheorem{corollary}[theorem]{Corollary}
\newtheorem{definition}[theorem]{Definition}

\newtheorem{remark}[theorem]{Remark}

\newcommand{\bpf}[1][Proof.]{\smallskip\noindent{\it #1} }
\newcommand{\qed}{\nolinebreak\mbox{\hspace{5 true pt}%
  \rule[-0.85 true pt]{3.9 true pt}{8.1 true pt}}}
\newcommand{\epf}{\qed \medskip}

\newtheorem{claim}[theorem]{Claim}
\newcommand{\cqed}{\nolinebreak\mbox{\hspace{5 true pt}%
  \rule[-0.85 true pt]{2.0 true pt}{8.1 true pt}}}
\newcommand{\bcpf}{\bpf[Proof of Claim.]}
\newcommand{\ecpf}{\cqed \medskip}

\begin{document}


\newenvironment{claimproof}[1]{\par\noindent\underline{Proof:}\space#1}{\leavevmode\unskip\penalty9999 \hbox{}\nobreak\hfill\quad\hbox{$\blacksquare$}}
\setcounter{tocdepth}{10}

\newcommand{\dedit}{\delta_{\mathrm{edit}}}
\newcommand{\Dedit}{\Delta_{\mathrm{edit}}}
\newcommand{\Done}{\hat{\Delta}_{1}}
\newcommand{\done}{\hat{\delta}_{1}}
\newcommand{\CP}{\mathcal{P}}
\newcommand{\PC}{\overline{\mathcal{P}}}

\newcommand{\join}{\nabla^{\bullet}}
\newcommand{\flip}{\nabla^{\bullet\bullet}}
\newcommand{\lambdamin}{\lambda_{\mathrm{min}}}
\newcommand{\lambdamax}{\lambda_{\mathrm{max}}}
\newcommand{\gammamax}{\gamma_{\mathrm{max}}}
\newcommand{\supp}{\mathrm{supp}}
\newcommand{\eps}{\varepsilon}
\newcommand{\OPT}{\mathrm{OPT}}
\newcommand{\Syone}{(Sym1)}
\newcommand{\Sytwo}{(Sym2)}
\newcommand{\Sone}{(Str1)}
\newcommand{\Stwo}{(Str2)}
\newcommand{\Stwoone}{(Str2.1)}
\newcommand{\Stwotwo}{(Str2.2)}
\newcommand{\xsupp}{(P1)}
\newcommand{\xbd}{(P2)}
\newcommand{\xone}{(P3)}
\newcommand{\xtwo}{(P4)}
\newcommand{\xcts}{(P5)}
\newcommand{\xH}{(P6)}

\title{Stability from graph symmetrisation arguments with applications to inducibility}
\author[1]{Hong Liu\footnote{Supported by UK Research and Innovation Future Leaders Fellowship MR/S016325/1 and
Institute for Basic Science Grant IBS-R029-C4}}
\author[2]{Oleg Pikhurko\footnote{Supported by 
ERC Advanced Grant 101020255 and
Leverhulme Research
Project Grant RPG-2018-424.}}
\author[3]{Maryam Sharifzadeh}
\author[4]{Katherine Staden\footnote{Supported by EPSRC Fellowship EP/V025953/1}}
\affil[1]{Extremal Combinatorics and Probability Group (ECOPRO), Institute for Basic Science (IBS), 55 Expo-ro, Yuseong-gu, Daejeon, 34126, South Korea}
\affil[2]{Mathematics Institute and DIMAP,
University of Warwick,
Coventry CV4 7AL, UK}
\affil[3]{Department of Mathematics and Mathematical Statistics, Ume\r{a} University, 901 87 Ume\r{a}, Sweden}
\affil[4]{School of Mathematics and Statistics, The Open University, Walton Hall, Milton Keynes MK7 6AA, UK}
\maketitle

\begin{abstract}
We present a sufficient condition for the stability property of extremal graph problems that can be solved via Zykov's symmetrisation. Our criterion is stated in terms of an analytic limit version of the problem. We show that, for example, it applies to the inducibility problem for an arbitrary complete bipartite graph $B$, which asks for the maximum number of induced copies of $B$ in an $n$-vertex graph,
and to the inducibility problem for $K_{2,1,1,1}$ and $K_{3,1,1}$, the only complete partite graphs on at most five vertices for which the problem was previously open.
\end{abstract}

\section{Introduction and notation}

The notion of symmetrisation in graphs was introduced by Zykov in~\cite{Zykov52}. In its most basic form, symmetrisation is the process of considering two non-adjacent vertices $x$ and $y$ in a graph $G$, and replacing $x$ by a clone of $y$, i.e.\ a vertex $y'$ whose neighbourhood is the same as that of $y$.
Zykov used symmetrisation to reprove Tur\'an's theorem~\cite{Turan41}, as follows.
Let $G$ be an $n$-vertex $K_r$-free graph with the maximum number of edges.
Whenever there are non-adjacent vertices $x,y$ with $d_G(x) \leq d_G(y)$, we symmetrise by replacing $x$ by a clone of $y$. The graph obtained in this way is still $K_r$-free and has at least as many edges as $G$, and one can do this so that the final graph is complete partite.
Standard convexity arguments imply that there are $r-1$ parts of almost equal size, recovering Tur\'an's theorem.
A variation of this approach was employed by Motzkin and Straus~\cite{MotzkinStraus65} also to reprove Tur\'an's theorem.

Suppose one seeks to maximise (or minimise) a graph parameter $\lambda$ such that there is always a way to symmetrise any given non-adjacent pair in a graph without decreasing $\lambda$.
Then it suffices to only consider `totally symmetrised' (that is, complete partite) graphs to determine the maximum value of $\lambda$.
Bollob\'as~\cite{Bollobas76} used symmetrisation to show that the parameter which counts any linear combination of cliques is symmetrisable, a special case of which provides a lower bound for the minimal number of cliques in a graph of given order and size.

In this paper, we are interested in more general graph parameters $\lambda$ which do not decrease upon symmetrisation, in a specific sense we describe below.
Like the example above, a symmetrisable $\lambda$ is maximised (not necessarily uniquely) by a complete partite graph.
Our main result gives a sufficient condition for stability for symmetrisable functions, namely that any graph which almost maximises $\lambda$ looks very much like a complete partite graph.
In fact we prove the quantitatively sharper property of \emph{perfect stability}, a strong form of stability which additionally implies an exact result.

\subsection{The statement of the main result}\label{sec-pre}

In order to define precisely what we mean by symmetrisable functions and perfect stability, we need to introduce some notation.
We write $G=(V,E)$ for a graph with vertex set $V$ and edge set $E$, and let $v(G) := |V|$ and $e(G) := |E|$.
Given $X \subseteq V$, we write 
$$
G[X] := (X,\{xy \in E: x,y \in X\})
$$ 
for the graph induced by $G$ on $X$, and $G-X := G[V(G)\setminus X]$, and also $G-x := G-\lbrace x \rbrace$.
Write $N_G(x) := \{y \in V: xy \in E\}$.

Fix a positive integer $k \geq 3$. Let $\C G$ be the family of all finite graphs up to isomorphism and let $\C G_n$ consist of graphs with $n$ vertices. Let $\CP_n \subseteq \C G_n$ be the family of complete partite graphs on $n$ vertices. Suppose we have a function $\gamma:\C G_k\to \I R$. For a graph $G=(V,E)$ with $v(G)\geq k$, define
\begin{equation}\label{eq:gamma}
\lambda(G):={n\choose k}^{-1} \sum_{X\in {V\choose k}} \gamma(G[X])
\end{equation}
where ${V\choose k}$ is the collection of $k$-element subsets of $V$.
Thus $\lambda(G)$ is the expected value of $\gamma(G[X])$ where $X$ is a random $k$-subset of $V$.
We may also work with 
$$
\Lambda(G):=\sum_{X\in {V\choose k}} \gamma(G[X])={n\choose k}\lambda(G),
$$ 
which may be more convenient in some calculations.  
For a vertex $x\in V(G)$, define
\begin{eqnarray*}
	\Lambda(G,x)&:=& \Lambda(G)-\Lambda(G-x)\ =\ \sum_{X\subseteq {V\choose k}: X\ni x} \gamma(G[X]),\\
	\lambda(G,x)&:=& {n-1\choose k-1}^{-1} \Lambda(G,x).
\end{eqnarray*}
Thus $\lambda(G,x)$ is the conditional expectation of $\gamma(G[X])$ where $X$ is a random $k$-subset of $V$ conditioned on containing~$x$.

Let $\lambda(n)$ be the maximum of $\lambda(G)$ over all $n$-vertex graphs $G$ and define 
$$
\lambdamax:=\lim_{n\to\infty} \lambda(n).
$$ One can easily show that the limit exists. Note that the minimisation problem reduces to a maximisation one just by negating $\gamma$, so we will always consider maximising $\lambda$ here.
We can now define what it means for $\lambda$ to be symmetrisable.

\begin{definition}[Symmetrisability]\label{symdef}
	A function $\lambda$ given by~(\ref{eq:gamma}) is \emph{symmetrisable} if for every $\e>0$ there is $n_0>0$ such that the following two properties hold for every graph $G=(V,E)$ of order $n\geq n_0$:
	\begin{enumerate}
		\item[\emph{\Syone}] There is a sequence of graphs $G_0,G_1,\ldots,G_m$ on $V$ such that $G_0=G$; $G_m$ is complete partite and for every $i\in [m]$  we have $\lambda(G_{i-1})\leq \lambda(G_i)$ and $|E(G_{i-1})\bigtriangleup E(G_i)|<\e {n\choose 2}$.
		\item[\emph{\Sytwo}] If $G-z$ is complete partite with partite sets $V_1,\ldots,V_t$, then there is a sequence of graphs $G_0,G_1,\ldots,G_m$ on $V(G)$ such that $G_0=G$; $G_i-z=G-z$; $\lambda(G_{i-1}) \leq \lambda(G_i)$; $|E(G_{i-1})\bigtriangleup E(G_i)|\leq \e (n-1)$ for all $i\in[m]$; and for each $j \in [t]$, either $N_{G_m}(z)\supseteq V_j$ or \mbox{$N_{G_m}(z)\cap V_j=\varnothing$}. 
	\end{enumerate}
\end{definition}

Here is an example of a symmetrisable parameter. For graphs $F,G$ with $v(F)\leq v(G)$, let $P(F,G)$ be the number of $v(F)$-subsets of $V(G)$ that induce a subgraph isomorphic to $F$. Let $p(F,G)=P(F,G)/{v(G)\choose v(F)}$ be the \emph{induced density} of $F$ in $G$. Let $\lambda(G):=\sum_{1 \leq i \leq k}a_ip(K_i,G)$ for $a_1,\ldots,a_k \in \mathbb{R}$. If we let $\gamma(F)=\sum_{1 \leq i \leq k} a_ip(K_i,F)$ for $F\in \C G_k$, then~\eqref{eq:gamma} holds. (Indeed, for $v(G)\geq k\geq i$, we have $p(K_i,G)=\sum_{F\in \C G_k} p(K_i,F)p(F,G)$ which implies the statement.) As mentioned above, Bollobas~\cite{Bollobas76} showed that $\lambda(n)$ is attained on a complete partite graph and his proof shows that every such $\lambda$ is in fact symmetrisable (for more details and examples, see Section~\ref{sec-app}).
In Section~\ref{sec-intro-app} we will see a generalisation of this parameter.

Secondly, we define perfect stability.
The \emph{edit} and \emph{normalised edit} distances between graphs $G$ and $H$ of the same order $n$ are given by
 \begin{eqnarray*}
 \Done(G,H)&:=&\min_{\sigma \in S(G,H)} |E(H)\bigtriangleup E(\sigma(G))|,\quad
 \done(G,H):= \frac{2}{n^2}\,\Done(G,H),
 \end{eqnarray*}
 where $S(G,H)$ is the set of bijections from $V(G)$ to $V(H)$.
(We also write $S(X) := S(X,X)$.)
We further define $\Done(G,\mathcal{H}):=\min_{H \in \mathcal{H}}\Done(G,H)$ for a family $\mathcal{H}$ of graphs of order $n$, and define $\done(G,\mathcal{H})$ analogously.

\begin{definition}[Perfect stability]
A graph parameter $\lambda$ is \emph{perfectly stable} if there exists $C>0$ such that for every graph $G$ of order $n\geq C$ there is a complete partite graph $H$ of order $n$
such that
$$\done(G,H)\leq C(\lambda(n)-\lambda(G)).$$
\end{definition}

We say that a sequence $\V x = (x_1,x_2,\ldots)$ with $x_1 \geq x_2 \geq \ldots \geq 0$ and $\sum_{i \geq 1}x_i \leq 1$ is a \emph{maximiser}
if there exists a sequence $(H_n)_n$ of complete partite graphs such that, as $n\to\infty$, we have $v(H_n) \rightarrow \infty$, $\lambda(H_n)\to \lambdamax$ and for every $i\geq 1$ the number of vertices in the $i$-th largest part of $H_n$ is $(x_i+o(1))v(H_n)$.
Let $\OPT=\OPT(\lambda)$ be the set of maximisers.

In Section~\ref{se:ff} we will show that if $\OPT$ is a finite set, then there is $\beta>0$ such that, for every $\V x\in\OPT$ and every $i\geq 0$ the entry $x_i$ is either 0 or at least $\beta$.

Observe that, if $\lambda$ is perfectly stable, then the only graphs on which $\lambda$ is maximised are complete partite. Perfect stability has already been proved in several contexts, most notably in Tur\'an-type problems; for example by F\"uredi~\cite{Furedi15}, Norin and Yepremyan~\cite{NorinYepremyan17,NorinYepremyan18}, Pikhurko, Slia\v{c}an and Tyros~\cite{PikhurkoSliacanTyros19}, and Roberts and Scott~\cite{RobertsScott18}.

\begin{definition}[Realisation $G_{n,\V x}$]\label{def:Gnx}
Given $n\in\mathbb{N}$ and $\V x = (x_1,x_2,\ldots)$ with $x_1 \geq x_2 \geq \ldots \geq 0$ and $x_0 := 1-\sum_{i \geq 1}x_i \geq 0$,
define a complete partite graph $G_{n,\V x}$ with vertex set $[n]$, parts $V_1,\ldots,V_m$ for some $m$ and a set $V_0$ of universal vertices, i.e., $|V_0|$ singleton parts, as follows.
If $x_0=0$, take a partition $[n]=V_1 \cup \ldots \cup V_m$ with $\left| |V_i|-x_i n\right| < 1$
and let $V_0=\varnothing$.
Otherwise, for all $i \geq 1$ with $x_i n \geq 2$, let $|V_i|=\lfloor x_i n \rfloor$ and let $V_0$ consist of the remaining vertices in $[n]$.

We say that $G_{n,\V x}$ is \emph{the ($n$-vertex) realisation} of $\V x$ and \emph{has $\PC$-structure} $V_0,\ldots,V_m$.
\end{definition}
If $H$ is a graph obtained by adding a new vertex $z$ to $G=G_{n, \V x}$, we say that $z$ is a \emph{clone} of $u \in V(G)$ if $u \in V_0$ and $N_H(z)=V(G)$, or if $u \notin V_0$ and $N_H(z)=N_G(u)$.
The following is one version of our main result, which is also stated as Theorem~\ref{thm-main}, in terms of limits.
Roughly speaking, it states that a symmetrisable function $\lambda$ is perfectly stable if it is `strict', meaning that it is sensitive to small alterations in a graph.

\begin{theorem}\label{thm-main-rough}
Let $\lambda$ be a symmetrisable function defined as above.
Suppose $|\OPT|<\infty$.
Suppose also that there exists $c>0$ such that the following hold for all large $n$ and maximisers $\V x=(x_1,x_2,\ldots) \in \OPT$, where $G=G_{n,\V x}$:
\begin{itemize}
\item[(i)] For all distinct $x,y \in V(G)$ we have $\lambda(G)-\lambda(G \oplus xy) \geq cn^{-2}$, where $G \oplus xy$ has vertex set $V(G)$ and edge set $E(G) \bigtriangleup \{\{x,y\}\}$.
\item[(ii)] 
If $G_v$ is obtained from $G$ by adding a new vertex $v$ which is complete or empty to each part of $G$
(where each $V_i$, $i \in [m]$ is a part and we have $|V_0|$ singleton parts)
then the minimum number of edits at $v$ needed to make $v$ a clone of some existing vertex of $G$ is at most $n(\lambda(G)-\lambda(G_v, v))/c$.
\end{itemize}
Then $\lambda$ is perfectly stable.
\end{theorem}

As mentioned, see Theorem~\ref{thm-main} for the `limit version' of this statement,
which concerns
$$
\lambda(\V x) := \lim_{n\to \infty}\lambda(G_{n,\V x}).
$$
One can easily show that this limit exists and that it does not depend on the choice of the part sizes $|V_i|$ in Definition~\ref{def:Gnx} (only on the ratios $x_i$).
The conditions in Theorem~\ref{thm-main-rough}~become a series of inequalities that must be verified for maximisers $\V x$, which are a finite collection of polynomial inequalities if the number of maximisers is finite and $x_0=0$,
since, for example, given $i,j$, the quantity $\lambda(G)-\lambda(G \oplus xy)$ is identical for all $x \in V_i$ and $y \in V_j$.
The value of the theorem is that, given the set of maximisers, the conditions are usually very easy to check, 
so in some sense the `combinatorial part' of the problem is solved.
It remains to determine the set of maximisers, amounting to a polynomial optimisation, which is unfortunately difficult in general.

\subsection{Applications to inducibility}\label{sec-intro-app}

A large class of problems where symmetrisation was sucessfully applied is the inducibility problem for complete partite graphs. The \emph{inducibility} problem for a graph $F$ is to determine $i(F,n):=\max\{P(F,G)\mid v(G)=n\}$, the maximum number of induced copies of $F$ that an order-$n$ graph $G$ can have. Note that $p(\overline{F},\overline{G})=p(F,G)$, 
where $\overline{G}$ denotes the complement of $G$,
so $i(\overline{F},n)=i(F,n)$. Also, consider
$$
i(F):= \lim_{n\to\infty} \frac{i(F,n)}{\binom{n}{v(F)}};
$$
the limit is known to exist and is in fact equivalent to the maximum density of induced copies of $F$ in a graphon $W$. Brown and Sidorenko~\cite[Proposition~1]{BrownSidorenko94} used symmetrisation to prove that if $F$ is
complete partite, then for every $n\in\I N$ at
least one
$i(F,n)$-extremal graph is complete partite.
Schelp and Thomason~\cite{SchelpThomason98}, also via symmetrisation, extended both the result of Brown and Sidorenko and a result of Bollob\'as~\cite{Bollobas76} by showing that the same conclusion holds (at least one graph attaining $\lambda(n)$ is complete partite) if the objective function is $\lambda(G)=\sum_F c_F\cdot p(F,G)$, where each $F$ is complete partite, including $K_t$ and $\overline{K}_t$, and $c_F$ is non-negative if $F$ is not a clique.
Their proof (which is essentially the same as that of Bollob\'as~\cite{Bollobas76}) implies that this parameter is symmetrisable (see Section~\ref{sec-app} for a proof).

\begin{lemma}[\cite{SchelpThomason98}]\label{lem-symm}
	The function  $\lambda(G):=\sum_{F} 
	c_F\cdot  p(F,G)$ is symmetrisable, where each $F$ is complete partite and $c_F\geq 0$ if $F$ is not a clique.\end{lemma}

In particular, Theorem~\ref{thm-main-rough} applies to the inducibility problem for complete partite graphs.
To the best of our knowledge, for every instance of this problem where the set of maximisers is known, we can prove perfect stability.

Pippenger and Golumbic~\cite{PippengerGolumbic75} determined $i(K_{s,t},n)$ for all $s,t$ with
$|s-t|\leq 1$, observing that the complete balanced bipartite graph is an
extremal graph.
Some of these results were independently reproved in~\cite{BollobasNaraTachibana86}.
Brown and Sidorenko~\cite{BrownSidorenko94} showed that $i(K_{s,t},n)$ with $st\geq 2$ is attained by a complete bipartite graph, and that if $\binom{t-s}{2}\leq s \leq t$ then the unique maximiser is $(\frac{1}{2},\frac{1}{2},0,\ldots)$. Perhaps surprisingly this does not mean that $K_{\lfloor n/2\rfloor,\lceil n/2\rceil}$ is optimal for $i(K_{s,t},n)$, and they show that if $3n = 4a^2+4$ for a large integer $a$, then $K_{n/2-a,n/2+a}$ is optimal for $K_{3,1}$.
We prove a corresponding stability result for complete bipartite graphs.

\begin{theorem}\label{thm-Kst}
Let $s,t \in \mathbb{N}$ with $st \geq 2$.
Then $p(K_{s,t},\cdot)$ is perfectly stable, $i(K_{s,t})=\binom{s+t}{s}M_{s,t}$ and there is a unique maximiser
$(\alpha,1-\alpha,0,0,0,\ldots)$,
where $\alpha \in [\frac{1}{2},1]$ maximises
$$
f_{s,t}(\alpha) := \alpha^s(1-\alpha)^t+\alpha^t(1-\alpha)^s
$$
and $M_{s,t}:=\max_{x\in[\frac{1}{2},1]} f_{s,t}(x)$ for $s\neq t$, 
and $M_{s,s}:=\frac{1}{2}\max_{x\in[\frac{1}{2},1]} f_{s,s}(x)$.
\end{theorem}

Bollob\'as, Egawa, Harris and Jin~\cite{BollobasEgawaHarrisJin95} studied the inducibility problem for complete equipartite graphs.
They showed that if the size $t$ of each part is not too small compared to the number $r$ of parts, then the complete balanced $r$-partite graph $T_r(n)$ is the \emph{unique} extremal graph for each large $n$. 
This strengthened an earlier work of
Brown and Sidorenko~\cite{BrownSidorenko94} which showed that $T_r(n)$ is an asymptotically extremal construction (without proving any uniqueness) -- that is, $(\frac{1}{r},\ldots,\frac{1}{r},0,\ldots)$ with $\frac{1}{r}$ repeated $r$ times is an element of $\OPT$.
We prove a corresponding stability result.

\begin{theorem}\label{thm-Krt}
Let $r,t \geq 2$ be integers and let $K_r(t)$ denote the complete $r$-partite graph with parts of size $t$.
Suppose that
$
t > 1+\log r
$
(denoting the natural logarithm by $\log$).
Then $p(K_r(t),\cdot)$ is perfectly stable, $i(K_r(t))=\frac{(tr)!}{r!t!^r r^{tr}}$, and the unique maximiser is $(\underbrace{\textstyle{\frac{1}{r},\ldots,\frac{1}{r}}}_{r},0,\ldots)$.
\end{theorem}
Interestingly, if the above lower bound on $t$ in terms of $r$ does not hold, then $(\frac{1}{r},\ldots,\frac{1}{r},0,\ldots) \notin \OPT$ (see~\cite{BrownSidorenko94}).

Finally, we obtain perfect stability for every previously unknown complete partite graph $F$ on $k \leq 5$ vertices. For this, note that trivially $K_k$ and $\overline{K}_k$ have unique maximisers $(0,0,\ldots)$, $(1,0,\ldots)$ respectively. If $F=K_{s,t}$ is bipartite, then Theorem~\ref{thm-Kst} implies that the unique maximiser $(\alpha,1-\alpha,0,\ldots)$ maximises
$
\alpha^s(1-\alpha)^t + \alpha^t(1-\alpha)^s
$.
Solving this, we see that $p(K_{s,t},\cdot)$ has unique maximiser $(\frac{1}{2},\frac{1}{2},0,\ldots)$ for all $s+t \leq 5$ apart from $\{s,t\}=\{4,1\}$, and here $p(K_{4,1},\cdot)$ has unique maximiser $(\frac{4}{5},\frac{1}{5},0\ldots)$.
Pikhurko, Slia\v{c}an and Tyros~\cite{PikhurkoSliacanTyros19} showed that $K_{2,1,1}$ is perfectly stable with unique maximiser $(\frac{1}{5},\ldots,\frac{1}{5},0,\ldots)$, and that
$K_{2,2,1}$ is perfectly stable with unique maximiser $(\frac{1}{3},\frac{1}{3},\frac{1}{3},0,\ldots)$ (we can also recover these results but do not provide proofs here).
The remaining $F$ are $K_{3,1,1}$ and $K_{2,1,1,1}$.
Flag algebra calculations of Even-Zohar and Linial~\cite{EvenZoharLinial15} give numerical upper bounds for these $i(F)$. Also, they provided lower bound constructions; these appear to match for both $K_{3,1,1}$ and $K_{2,1,1,1}$.
They speculated that their lower bound constructions are tight in both cases.
We confirm this and prove perfect stability for these $F$. (After this paper was submitted, Liu, Mubayi and Reiher~\cite[Theorem 1.13]{LiuMubayiReiher23} determined the value of $i(K_t^-)$ for every $t$, where $K_t^-=K_{2,1,\dots,1}$ is the complete graph of order $t$ minus one edge.)

\begin{theorem}\label{thm-K2111}
$p(K_{2,1,1,1},\cdot)$ is perfectly stable, $i(K_{2,1,1,1}) = \frac{525}{1024}$, and the unique maximiser is $(\frac{1}{8},\ldots,\frac{1}{8},0,\ldots)$.
\end{theorem}

\begin{theorem}\label{thm-K311}
$p(K_{3,1,1},\cdot)$ is perfectly stable, $i(K_{3,1,1}) = \frac{216}{625}$, and the unique maximiser is $(\frac{3}{5},0,\ldots)$.
\end{theorem}

The latter is particularly interesting since the extremal graph contains a clique part: it is a clique with a clique of proportion $3/5$ removed. This demonstrates that allowing maximisers $\V x$ with $x_0=1-\sum_{i \geq 1}x_i>0$ in our theory -- which complicates matters somewhat -- is essential in giving a full picture.

We remark that the case $\lambda(\cdot)=-p(\overline{K_3},\cdot) - p(K_3,\cdot)$ (which is not a function as in Lemma~\ref{lem-symm}) is given by a classical theorem of Goodman, who determined this value exactly.
Here, asymptotically extremal graphs are those for which all but $o(n)$ vertices have degree $\frac{n}{2}+o(n)$ (including many graphs which are not complete partite).
(Note that $p(\overline{K_3},\cdot) + p(K_3,\cdot)$ is trivially maximised by the complete and empty graphs.)
It remains a major open problem to determine $\lambdamax$ for $\lambda(\cdot)=-p(\overline{K_4},\cdot)-p(K_4,\cdot)$.

Pikhurko, Slia\v{c}an and Tyros~\cite{PikhurkoSliacanTyros19}
were able to prove perfect stability for $i(F,n)$ for several small graphs $F$ via flag algebra calculations. The graphs they considered were 
$C_4=K_{2,2}$, $K_{2,1,1}$, $K_{3,2}$, $K_{2,2,1}$, as well as the non-complete partite graphs $P_3 \cup K_2$, the ``$Y$" graph and the paw graph which we do not define.
Their results extend inducibility results obtained in~\cite{BrownSidorenko94},~\cite{PippengerGolumbic75}, and by Hirst in~\cite{Hirst14}.
Our Theorem~\ref{thm-Kst} in particular reproves the cases $K_{2,2}$ and $K_{3,2}$ from~\cite{PikhurkoSliacanTyros19}.

Before stating the limit version of our main theorem in Section~\ref{sec-strict}, we give here an illustration of it in the case $F=C_4$.
(Perfect stability was already proved here in~\cite{PikhurkoSliacanTyros19}.)
It is easy to see that $\OPT$ consists only of the unique vector $(\frac{1}{2},\frac{1}{2},0,\ldots)$ with $\lambdamax=\frac38$. Thus in order to apply our criterion we have to check that, starting with $K_{\lfloor n/2\rfloor,\lceil n/2\rceil}$ the following two properties hold:  (i) if we add an edge into a part or remove an edge across then we decrease the number of induced copies of $C_4$ by $\Omega(n^2)$; (ii) if we add a new vertex $v$ which is either isolated or connected to every other vertex, the number of induced copies of $C_4$ containing $v$ is at most $(1-\Omega(1)) \,\frac38\binom{n}{3}$. Both properties trivially hold so the inducibility problem for $C_4$ is indeed perfectly stable by Theorem~\ref{thm-main-rough}.

The following conjecture seems plausible.

\begin{conjecture}
The inducibility problem for $F$ is perfectly stable for every complete partite~$F$.
\end{conjecture}

However it is not the case that every problem with $\lambda = \sum_F c_F \cdot p(F,\cdot)$ is perfectly stable, where each $F$ is complete partite, and $c_F \geq 0$ if $F$ is not a clique.
Indeed, if $k \geq 3$ and the sum is over all complete partite $F$ on $k$ vertices, and each $c_F = 1$, then every $k$-vertex subset of every complete partite graph contributes (the maximum value of) $1$ to $\Lambda$, so $\OPT$ is the set of \emph{all} $\V x$ with $x_1 \geq x_2 \geq \ldots \geq 0$ and $\sum_{i \geq 1}x_i \leq 1$.
Let us show that $\lambda$ is not perfectly stable.
Indeed, if it is, there is $C$ such that for every graph $G$ of order $n \geq C$, there is a complete partite $H$ such that $\done(G,H) \leq C(\lambda(n)-\lambda(G))$.
Choose $1/n \ll c \ll 1/C$.
Starting with $K_n$, remove every edge with both endpoints inside a set $A$ of size $5cn$ and add into $A$ a blow-up of $C_5$ with each part $A_1,\ldots,A_5$ of size $cn$, to obtain an $n$-vertex graph $G$.
Then $\done(G,H) = \Omega(c^2)$ for every complete partite $H$,
but $\lambda(n)-\lambda(G)=1-\lambda(G) = O(c^3)$.
Indeed, a subset of $G$ is not complete partite only if it contains at least three vertices in $A$. So the fraction of subsets inducing a non-complete partite graph is $O(c^3)$.
This is a contradiction.

Finally, it would be remiss not to remark on the inducibility problem for non-complete partite graphs, for which the present paper does not apply, and which is in general wide open (see~\cite{EvenZoharLinial15} for a list of known results of order up to $5$). 
The outstanding open problem in the area is determining $i(P_4)$, the smallest unsolved case, for which there is not even a conjectured value.
Hatami, Hirst and Norin proved that extremal graphs of large blow-ups are essentially blow-ups themselves~\cite{HatamiHirstNorin14}.
Graphs with more interesting structure appear as extremal graphs for other $F$.
An important longstanding conjecture of Pippenger and Golumbic~\cite{PippengerGolumbic75} is that $i(C_k) = k!/(k^k-k)$ for $k \geq 5$, attained by the \emph{iterated blow-up} of $C_k$.
Balogh, Hu, Lidick\'y and Pfender~\cite{BaloghHuLidickyPfender16} proved this conjecture for $k=5$: they obtained an exact result for $\lambda(\cdot)=p(C_5,\cdot)$ and showed that if $n$ is a power of $5$ then the unique graph attaining $i(C_5,n)$ is an iterated blow-up of a $5$-cycle.
There has recently been progress on the general conjecture, see~\cite{HefetzTyomkyn18,KralNorinVolec19}.
Yuster~\cite{Yuster19} and independently Fox, Huang and Lee~\cite{FoxHuangLee} proved that for almost all graphs $F$, the extremal graph is the iterated blow-up of $F$.
Fox, Sauermann and Wei~\cite{FoxSauermannWei} considered graphs $H$ obtained by removing a small number of vertices from a random Cayley graph $\tilde{H}$ of an abelian group, showing that here the extremal graph is the iterated blow-up of $\tilde{H}$ (not of $H$).
Liu, Mubayi and Reiher~\cite{LiuMubayiReiher23} began a systematic study of the feasible region of induced graphs; that is, the set of points $(x,y)$ in the unit square for which there is a graph of edge density approaching $x$ with $F$-density approaching $y$. The inducibility problem seeks the maximum $y$-value of such a point.

The directed analogue of the inducibility problem is also actively studied, e.g.~for stars~\cite{Huang14,HuMaNorinWu}, paths~\cite{ChoiLidickyPfender20}
and $4$-vertex graphs~\cite{BozykGrzesikKielak,BurkeLidickyPfenderPhillips,HuLidickyPfenderVolec}.

\subsection{Structure of the paper}


The rest of the paper is organised as follows.
In Section~\ref{sec-space} we introduce the partite limit space corresponding to the collection of limits of complete partite graphs
which we will need to prove our main result.
In Section~\ref{sec-strict} we define the notion of strictness in terms of elements of this space and give a limit version of our main result, Theorem~\ref{thm-main}.
The main result of Section~\ref{se:ff} is that when $\OPT$ is finite, all part ratios of extremal graphs are bounded away from $0$.
We prove Theorem~\ref{thm-main} in Section~\ref{sec-main}.
We present some applications of Theorem~\ref{thm-main} to the inducibility problem (Theorems~\ref{thm-Kst}--\ref{thm-K311}) in Section~\ref{sec-app}.
Section~\ref{sec-conclude} contains some concluding remarks.

We denote by $\mathbb N:=\{1,2,\dots\}$ and $\mathbb N_0:=\{0,1,\dots\}$ the sets of respectively positive and non-negative integers.


\section{The partite limit space}\label{sec-space}

We will work in a space $\PC$, the \emph{partite limit space}, which is in some sense the completion of the set of complete partite graphs.
The aim of this section is to define $\PC$ and a metric $\dedit$ on this set, which will essentially generalise edit distance in graphs. We prove that this yields a compact metric space upon which $\lambda$ can be extended continuously (Lemma~\ref{lm:compact}).
Thus the set $\OPT$ of maximisers of $\lambda$ in $\PC$ is non-empty.
We define
$$
\PC := \left\lbrace \V x = (x_1,x_2,\ldots) : x_1 \geq x_2 \geq \ldots \geq 0 \text{ and } \sum_{i \geq 1}x_i\leq 1\right\rbrace.
$$
As usual, $\supp(\V x) := \lbrace i \geq 1: x_i>0\rbrace$, and we also define $\supp^*(\V x) := \supp(\V x) \cup \lbrace 0 \rbrace$ if $\sum_{i\geq 1} x_i<1$, and $\supp^*(\V x) := \supp(\V x)$ otherwise.
For $\beta>0$ we write 
$$
\PC_\beta := \{ \V x \in \PC: x_i \geq \beta \ \forall \ i \in \supp^*(\V x)\}.
$$
Write ${\bf 0} := (0,0,\ldots)$.
Given $\V x, \V x_n \in \PC$, we will always write $\V x = (x_1,x_2,\ldots)$ and $\V x_n = (x_{n,1},x_{n,2},\ldots)$ and correspondingly $x_0:=1-\sum_{i \geq 1}x_i$ and $x_{n,0}:=1-\sum_{i \geq 1}x_{n,i}$.
A complete partite graph $G=K(V_1,\ldots,V_m)$ on vertex set $[n]$ with $|V_1|\geq \ldots\geq |V_m|$ corresponds to the 
vector 
 $$
 \V  x_G:=(|V_1|/n,\ldots,|V_m|/n,0,\ldots).
 $$ 
 We write $\CP$ for the set of those elements $\V x$ of $\PC$ with finitely many non-zero entries all of which are rational, thus corresponding to the set of complete partite graphs.
 Somewhat conversely, we have the construction $G_{n,\V x}$ from Definition~\ref{def:Gnx}.
For example, we have $G_{n,{\bf 0}} \cong K_n$, $G_{n,(1,0,\ldots)} \cong \overline{K}_n$ and (assuming $n=2\ell$ is even)  $G_{n,(\frac{1}{2},\frac{1}{2},0,\ldots)} \cong K_{\ell,\ell}$,
but we cannot take, say, any $K_{a,b,1}$ for $G_{n,(x,1-x,0,\ldots)}$.

\subsection{The measure-theoretic and graphon perspectives}\label{se:graphon}

For each $\V x \in \PC$, one can define a probability measure $\mu_{\V x}$ on $\mathbb{N}_0$
by setting $\mu_{\V x}(\{i\})=x_i$ and then let
$$
\mathcal{M} := \left\{ \mu_{\V x}: \V x \in \PC\right\}.
$$

It is very natural to define the corresponding collection of ``complete partite" graphons (which will be used in Section~\ref{se:ff}).
A \emph{graphon} is a quadruple $Q=(\Omega,\mathcal{B},\mu,W)$,
where $(\Omega,\mathcal{B},\mu)$ is a standard probability space and $W: \Omega \times \Omega \to [0,1]$ is a symmetric measurable function.
For every graph $G$, we define the corresponding graphon $Q_G=(V,2^V,\mu,A_G)$
where $\mu$ is the uniform measure on the finite set $V$ and $A_G: V\times V \to \{0,1\}$ is the adjacency function of $G$.
For a graph $F$ on $[k]$ we write
$$
p(F,Q) := \frac{k!}{|\rm{aut}(F)|}\int_{\Omega^k}\prod_{ij \in E(F)}W(x_i,x_j)\prod_{ij \in E(\overline{F})}(1-W(x_i,x_j))d\mu(x_1)\ldots d\mu(x_k)
$$
where ${\rm{aut}}(F)$ is the group of automorphisms of $F$.
In the literature one usually encounters $t_{\rm ind}(F,Q)$ which is the above without the normalisation factor.
Two graphons $Q,Q'$ are \emph{equivalent} or \emph{weakly isomorphic} if $p(F,Q)=p(F,Q')$ for every graph $F$.
A sequence of graphons $(Q_n: n \in \mathbb{N})$ is said to \emph{converge} to a graphon $Q$ if $\lim_{n\to\infty}p(F,Q_n)=p(F,Q)$ for every graph $F$.
A \emph{$Q$-random graph} of order $k$ is obtained by sampling $k$ random points $v_1,\ldots,v_k \in (\Omega,\mu)$ uniformly and independently, and adding each edge $x_ix_j$ with probability $W(x_i,x_j)$.

Now let $Q_{\V x} := (\mathbb{N}_0, 2^{\mathbb{N}_0},\mu_{\V x}, K)$
where $K(i,j):=0$ if $i=j\geq 1$ and $K(i,j):=1$ otherwise,
i.e.\ if $i \neq j$ or $i=j=0$. Then define
$$
\mathcal{Q} := \left\{Q_{\V x}: \V x \in \PC\right\}.
$$
There are various characterisations of weak isomorphism (see~\cite[Theorem~13.10]{Lovasz12}).
All we will need is the easy fact that for distinct $\V x, \V y \in \PC$, their graphons $Q_{\V x}, Q_{\V y}$ are not weakly isomorphic.
Indeed, if $i\geq 1$ is the minimum integer with $x_i\not=y_i$, say
$x_i>y_i$, then it is not hard to see directly that the edgeless graph
of sufficiently large order $n$ has strictly larger density in $\V x$
than in $\V y$.

The spaces $\PC$, $\mathcal{M}$ and $\mathcal{Q}$ are equivalent and one can take any of these perspectives, but in this paper we mainly work with $\PC$ (and briefly use $\mathcal{Q}$ in Section~\ref{se:ff}).
The space $\mathcal{Q}$ was used in~\cite{BennettDudekLidickyPikhurko20} by Bennett, Dudek, Lidi\'cky and Pikhurko who determined the minimum $C_5$-density in graphs of edge density $\frac{k-1}{k}$ for integers $k$. They used $\mathcal{Q}$ to prove a corresponding stability result.
Therefore we hope that the theory concerning $\PC$ (and, by extension, $\mathcal{M}$ and $\mathcal{Q}$) developed in this section may be useful for other extremal problems where the extremal graphs are complete partite.

\subsection{The edit metric}

We would like to define a metric on $\PC$ which will correspond to the edit distance between graphs.
First we define edit distance between two graphs of possibly different orders, often called the \emph{fractional edit distance}. Given a graph $G$, let $G^{(n)}$ be an $n$-vertex almost uniform blow-up of $G$, that is, we replace each vertex $x \in V(G)$ by an independent set $I_x$, where each $\left| |I_x|-|I_y|\right| \leq 1$, and $\sum_{x \in V(G)}|I_x|=n$, and add every edge between $I_x$ and $I_y$ whenever $xy \in E(G)$. Then let
$$
\dedit(G,H) := \lim_{n\to\infty}\done(G^{(n)},H^{(n)}).
$$
It is easy to see that the limit exists; in fact, its value can be
computed via a linear program with $v(G)\times v(H)$ variables that
considers all fractional overlays  between the vertex sets of $G$ and
$H$, c.f.\ e.g. \cite[Equation (3)]{Pikhurko10}.
We also define for a family $\mathcal{H}$ of graphs $\dedit(G,\mathcal{H}) := \lim_{n \to \infty}\done(G^{(n)},\{H^{(n)}: H \in \mathcal{H}\})$.
We define the distance between $\V x,\V y\in\PC$ to be
 $$
 \dedit(\V x, \V y) := \lim_{n\to\infty}\done(G_{n,\V x},G_{n,\V y}).
 $$
 For a graph $G$, define also $\dedit(\V x, G) := \lim_{n\to\infty}\done(G_{n,\V x}, G^{(n)})$ and, for a family $\mathcal{H}$ of graphs, $\dedit(\V x, \mathcal{H}))$ in the obvious way. Again, the existence of the limit in these definitions is easy to establish.
Note that the normalisation factor $\frac{2}{n^2}$ in the `usual' edit distance $\done$
is motivated by vertices of $G$ corresponding to independent sets of relative size $\frac{1}{n}$.
The distances $\done$ and $\dedit$ are not the same even for graphs of the same order, due to rounding; see examples of Matsliah (see Appendix~B in \cite{GoldreichKrivelevichNewmanRozenberg08}) and Pikhurko~\cite{Pikhurko10}.
The following lemma implies that we are free to interchange $\dedit$ and $\done$ in matters of convergence, and that with respect to $\dedit$ we are free to interchange $H$ and $\V x_H$ when $H$ is complete partite.
\begin{lemma}\label{lm:edit}
We have the following.
\begin{itemize}
\item[(i)] $\dedit(G,H) \leq \done(G,H) \leq 3\dedit(G,H)$ for graphs $G,H$ with the same order.
\item[(ii)] $\dedit(H, \V x_H)=0$ and $\dedit(\V x, H) = \dedit(\V x, \V x_H)$ and $\dedit(G,H)=\dedit(\V x_G, \V x_H)$ for all $\V x \in \PC$ and complete partite graphs $G,H$.
\item[(iii)] $\dedit$ satisfies the triangle inequality on $\PC$.
\end{itemize}
\end{lemma}

\bpf
The non-trivial inequality of part~(i) was proved in \cite[Lemma~14]{Pikhurko10}.
For~(ii), let $H$ have $h$ vertices. Then $\V x_{H^{(nh)}}=\V x_H$ and $G_{nh,\V x_H} = (G_{h,\V x_H})^{(n)}=H^{(nh)}$ for any integer $n$. Since any subsequence of $\left(\done(H^{(m)},G_{m,\V x_H})\right)_m$ converges to $\dedit(H,\V x_H)$, we have
$$
\dedit(H,\V x_H)= \lim_{n\to\infty}\done(H^{(nh)},G_{nh,\V x_H})=\lim_{n\to\infty}\done(H^{(nh)},H^{(nh)})=0.
$$ 
The remaining parts of~(ii) now follow from~(iii)
which is immediate since $\done$ satisfies the triangle inequality on the set of graphs of the same given order.
\epf

This notion of edit distance is very natural, yet rather unwieldly to work with. The following easy facts concerning it will be useful.
Recall first that $x_0$ is not an entry in $\V x = (x_1,x_2,\ldots)$, so e.g.~$\| \V x\|_1 = \sum_{i \geq 1}|x_i| = 1-x_0$.

\begin{proposition}\label{editprops}
For all $\V x, \V y \in \PC$, we have that 
\begin{itemize}
\item[(i)] $\dedit(\V x, \V y) \leq 2\| \V x - \V y \|_1$.
\item[(ii)] $\dedit(\V x, \bf{0}) = \| \V {x} \|_2^2$.
\item[(iii)] $\dedit(\V x, (x_1,\ldots,x_M,0,\ldots)) \leq \sum_{i>M}x_i^2$ for all $M \geq 1$.
\end{itemize}
\end{proposition}

\bpf
For (i), consider large $n \in \mathbb{N}$ and $G_{n,\V x}, G_{n,\V y}$ with $\PC$-structures $V_0,V_1,\ldots,V_m$ and $U_0,U_1,\ldots,U_\ell$ respectively, where without loss of generality $\ell \leq m$.
For convenience let $U_{\ell+1}=\ldots = U_m = \varnothing$.
Let $\sigma \in S([n])$ be a permutation (and recall that $V(G_{n,\V x})=V(G_{n, \V y})=[n]$).
For all $0 \leq i,j \leq m$,
let $X_{ij}=\sigma(V_i)\cap U_j$.
A pair of vertices is included in the symmetric difference $E(G_{n, \V x}) \bigtriangleup E(\sigma(G_{n,\V y}))$ if and only if either it 
lies in $V_i$ for some $i\ge 1$ but not in $X_{ij}$ for any $j\in [m]$,
or lies in $U_j$ for some $j\ge 1$ but not in $X_{ij}$ for any $i\in [m]$. Thus
\begin{equation}\label{eq:edit}
E(G_{n, \V x}) \bigtriangleup E(\sigma(G_{n,\V y}))=\sum_{i \in [m]}\left(\binom{|V_i|}{2}-\sum_{j \in [m]}\binom{|X_{ij}|}{2}\right)+\sum_{j \in [m]}\left(\binom{|U_j|}{2}-\sum_{i \in [m]}\binom{|X_{ij}|}{2}\right).
\end{equation}
Take $\sigma \in S([n])$ so that for all $i \geq 0$, $\sigma(V_i) \subseteq U_i$ whenever $|V_i| \leq |U_i|$, and $\sigma(V_j) \supseteq U_j$ whenever $|U_j|\leq |V_j|$ (and $\sigma$ is otherwise arbitrary).
If $|U_k|<|V_k|$ then $|X_{kk}|=|U_k|$ and $|X_{ik}|=0$ for all $i \neq k$ and thus
\begin{align*}
&\binom{|V_k|}{2} - \sum_{j \in [m]}\binom{|X_{kj}|}{2} + \binom{|U_k|}{2} - \sum_{j \in [m]}\binom{|X_{jk}|}{2}= \binom{|V_k|}{2} - \sum_{j \in [m]}\binom{|X_{jk}|}{2}\\
&\leq \binom{|V_k|}{2} - \binom{|X_{kk}|}{2} \leq \frac{1}{2}(|V_k|^2-|U_k|^2)+|U_k|.
\end{align*}
Then
\begin{align*}
\Done(G_{n,\V x}, G_{n,\V y}) &\leq \sum_{k \in [m]}\left(\frac{1}{2}\left|\, |V_k|^2-|U_k|^2\, \right|+\min\{|U_k|,|V_k|\}\right) \leq n\sum_{k \in [m]}\left(|\, |V_k|-|U_k|\, |+y_k+o(1)\right)\\
&\leq n^2\| \V x - \V y\|_1 + O(n).
\end{align*}
So $\dedit(\V x, \V y) \leq 2\|\V x - \V y\|_1$, as required.

Parts~(ii) and~(iii) are clear.
\epf

Note however that convergence in $\ell_1$ does not give the same topology as pointwise convergence, by considering for each $n \in \mathbb{N}$ the sequence $\V x_n$ given by $x_{n,i}=1/n$ for all $i \in [n]$ and $x_{n,i}=0$ otherwise.
We have that $\| \V x_n \|_1=1$ for all $n$, while $\V x_n$ clearly converges pointwise to $\bf{0}$ and by Proposition~\ref{editprops}(ii), we see that $\dedit(\V x_n, {\bf 0}) =\frac{1}{n}\rightarrow 0$ as $n\rightarrow \infty$. On the other hand, convergence in $\dedit$ is equivalent to pointwise convergence, as we show in the next lemma.

\begin{proposition}\label{lem-edit-pw}
	In the space $\PC$, convergence in edit distance is equivalent to pointwise convergence.
	That is, whenever $(\V x_n)_n$ is a sequence in $\PC$ and $\V x \in \PC$, we have that $\lim_{n \rightarrow \infty}\dedit(\V x_n,\V x)=0$ if and only if for all $i \in \mathbb{N}$ we have that $\lim_{n\rightarrow \infty}|x_{n,i}-x_i|=0$.
\end{proposition}
\bpf 
Let $(\V x_n)_n$ be a sequence in $\PC$ and let $\V x \in \PC$. Fix an arbitrary $\eps>0$. 

Suppose first that $\V x_n \rightarrow \V x$ pointwise.
We need to show that $\dedit(\V x_n, \V x) < \eps$ for sufficiently large $n$. Since $\sum_{i \geq 1}x_i \leq 1$ and $x_1\geq x_2\ge\ldots \geq 0$, there exists an integer $M>0$ such that $\sum_{i \geq M}x_i < \eps/8$, in particular, $x_M < \eps/8$. As $\V x_n \rightarrow \V x$ pointwise, there exists $n_0$ such that, for all $i \leq M$ and for all integers $n \geq n_0$, we have that $|x_{n,i}-x_i|<\eps/(8M)$.
In particular, since $x_{n,j}$ is non-increasing with $j$, we have for all integers $n \geq n_0$ and $j \geq M$ that $x_{n,j}<\eps/4$.
Let $\V y := (x_1,\ldots,x_M,0,\ldots)$ and, for each $n \in \mathbb{N}$, define $\V y_n := (x_{n,1},\ldots,x_{n,M},0,\ldots)$.
Let $n \geq n_0$ be an integer. Then by Proposition~\ref{editprops}(i),
$
\dedit(\V x, \V y) \leq 2\|\V x - \V y\|_1 = 2\sum_{i > M}x_i < \frac{\eps}{4}
$.
Similarly by Proposition~\ref{editprops}(ii) and~(iii),
$$
\dedit(\V x_n, \V y_n) \leq \dedit((x_{n,M+1},x_{n,M+2},\ldots), \mathbf{0}) = \sum_{i > M}x_{n,i}^2 \leq \sup_{i > M}x_{n,i} \cdot \sum_{i > M}x_{n,i} \leq x_{n,M+1} \leq \frac{\eps}{4}.
$$
But also
$
\dedit(\V y, \V y_n) \leq 2\|\V y - \V y_n\|_1 = 2\sum_{i \leq M}|x_{n,i}-x_i| < \frac{\eps}{4}
$.
By Lemma~\ref{lm:edit}(iii), $\dedit$ defined on $\PC$ satisfies the triangle inequality.
Thus we have $\dedit(\V x_n, \V x) \leq \eps$ whenever $n \geq n_0$.
Thus $\V x_n \rightarrow \V x$ in edit distance, as required.

Conversely, suppose now that $(\V x_n)_n$ converges to $\V x$ in edit distance $\dedit$. Let $i \geq 1$. We need to show that there exists $n_0>0$ such that for all $n>n_0$ we have $|x_{n,i}-x_i|\leq \e $.
Now, there exists $n_0>0$ such that for all $n>n_0$, there is a permutation $\sigma:[n]\rightarrow [n]$ such that
\begin{equation*}
\Done(G_{n,\V x_n},G_{n,\V x}) = |E(G_{n,\V x_n})\bigtriangleup E(\sigma(G_{n,\V x}))|\leq (\eps n/12)^2.
\end{equation*}
Let $n>n_0$. For $A\subseteq [n]$, denote by $\sigma(A)$ and $\sigma^{-1}(A)$ the image and pre-image of $A$ respectively. By definition $G_{n,\V x_n}$ has a vertex partition $V_{n,0} \cup V_{n,1} \cup \ldots \cup V_{n,m}$, where $V_{n,0}$ is a clique, $V_{n,i}$ is an independent set for all $i \in [m]$, and $G_{n,\V x_n}$ 
is complete between every distinct $V_{n,i}$ and $V_{n,j}$.
Define $V_0 \cup V_1 \cup \ldots \cup V_\ell$ analogously for $G_{n,\V x}$. So
\begin{equation}\label{Vdef}
|V_{n,i}| = x_{n,i}n + O(1) \quad\text{and}\quad |V_i| = x_i n + O(1) \quad\text{ for all } i \geq 0.
\end{equation}
Choose an ordering of the vertices of $G_{n,\V x_n}$ so that a vertex $u \in V_{n,i}$ comes before a vertex $v \in V_{n,j}$ if $1 \leq i < j$; or if $i \neq 0$ and $j=0$.
Choose an analogous ordering for $V(G_{n,\V x})$.
Note the following trivial equality:
\begin{eqnarray}\label{eq-imgpreimg}
|\sigma(V_{j})\bigtriangleup V_{n,i}|=|V_{j}\bigtriangleup \sigma^{-1}(V_{n,i})|\quad\text{ for all }\quad i,j\in\mathbb{N}_0.
\end{eqnarray}
We first show that for each vertex part $V_{n,i}$ which is not too small, there is a unique part $V_{j_i}$ such that $\sigma$ maps most of $V_{n,i}$ to $V_{j_i}$.
Given $i \in \{0,1,\ldots,m\}$ and $j \in \{0,1,\ldots,\ell\}$, we say that $i$ is \emph{$j$-good} if 
$
|\sigma(V_{j}) \bigtriangleup V_{n,i}| < \eps n/4
$.

\begin{claim}\label{cl-fullset}
Let $A := \lbrace i \in [m] : |V_{n,i}| \geq \eps n/2 \rbrace$.
Then there exists $B \subseteq [\ell]$ with $|A|=|B|$ and a bijection $\mu : A \rightarrow B$ such that, for every $i \in A$, we have that $i$ is $j$-good if and only if $j=\mu(i)$.
\end{claim}
\bcpf
Let $i \in A$.
Note first that $i$ is not $0$-good. Indeed, this follows from
$$
(\eps n/12)^2 \geq \Done(G_{n,\V x_n},G_{n,\V x}) \geq \binom{|\sigma(V_0) \cap V_{n,i}|}{2}.
$$
 So
$
	\sum_{j\in\mathbb{N}}|\sigma(V_j)\cap V_{n,i}|=|V_{n,i}|- |\sigma(V_0) \cap V_{n,i}| > \eps n/4
$.
Suppose now that $i$ is not $j$-good for any $j \in [\ell]$. Then since $G_{n,\V x}[V_j]$ and $G_{n,\V x_n}[V_{n,i}]$ are empty graphs, and both $G_{n,\V x}$ and $G_{n,\V x_n}$ are complete partite graphs, 
$$
|E(G_{n,\V x_n})\bigtriangleup E(\sigma(G_{n,\V x}))|\geq \sum_{j\in\mathbb{N}}|\sigma(V_j)\bigtriangleup V_{n,i}|\cdot |\sigma(V_j)\cap V_{n,i}|\geq \eps n/4 \cdot \sum_{j\in \mathbb{N}}|\sigma(V_j)\cap V_{n,i}|> (\eps n/4)^2,
$$
a contradiction. Thus there is some $j_i \in \mathbb{N}$ for which $i$ is $j_i$-good.
We claim that we can set $\mu(i) := j_i$ and $B := \lbrace \mu(i) : i \in A\rbrace$.
We first show that this is well-defined, i.e.\ $j_i$ is unique. Fix an arbitrary $j'\in [\ell]\setminus \{j_i\}$. 
Since $\sigma$ is a permutation, $\sigma(V_{j'})\cap \sigma(V_{j_i})=\varnothing$, and therefore
\begin{eqnarray*}
	|\sigma(V_{j'})\bigtriangleup V_{n,i}|\geq |\sigma(V_{j_i})\cap V_{n,i}|\geq |V_{n,i}|-|\sigma(V_{j_i})\bigtriangleup V_{n,i}|>\eps n/4,
\end{eqnarray*}
i.e.\ $i$ is not $j'$-good. It remains to show that $\mu$ is injective, i.e.\ that if $i' \in A\setminus \lbrace i \rbrace$, we have that $i'$ is not $j_i$-good.
By~\eqref{eq-imgpreimg}, it suffices to show that $|V_{j_i}\bigtriangleup \sigma^{-1}(V_{n,{i'}})|\geq \eps n/4$. Since $\sigma^{-1}$ is a permutation, $\sigma^{-1}(V_{n,i'})\cap \sigma^{-1}(V_{n,i})=\varnothing$, and therefore
$
	|V_{j_i}\bigtriangleup \sigma^{-1}(V_{n,i'})|\geq |V_{j_i}\cap \sigma^{-1}(V_{n,i})|>\eps n/4
$
as desired, where the last inequality follows from $i$ being $j_i$-good and~\eqref{eq-imgpreimg}. This completes the proof of the claim.
\ecpf

We are now ready to prove the desired conclusion that for all $i\in\I N$, $|x_{n,i}-x_i|\leq \e$. Suppose this is not true, and let $k$ be the smallest integer $i$ such that $|x_{n,i}-x_i|> \eps$. 
Assume that $x_{n,k}>x_k+\eps$ (the other case can be handled similarly). In particular, recalling~(\ref{Vdef}), $|V_{n,k}|\geq \eps n/2$, and so $[k]\subseteq A$.
Since $x_1 \geq x_2 \geq \ldots$ and $x_{n,1} \geq x_{n,2} \geq \ldots$, we have for all $1 \leq i \leq k \leq i'$ that, neglecting $O(1/n)$ error terms,
$x_{n,i} \geq x_{n,k}>x_{k}+\eps\geq  x_{i'}+\eps$, so
$
|V_{n,i}| \geq |V_{i'}|+\eps n/2
$.
Thus, for all positive integers $i \leq k$, we must have $\mu(i) < k$.
In other words, $\mu([k])\subseteq [k-1]$, which contradicts $\mu$ being a bijection.
This completes the proof of the lemma.
\epf 
\begin{remark} Lemma~3.8 in \cite{BennettDudekLidickyPikhurko20} proves that
if $\V x_n, \V x \in \PC$ are such that $\V x_n \to \V x$ pointwise, then the corresponding graphons $Q_{\V x_n}$ converge to $Q_{\V x}$; that is, all the $p(F,Q_{\V x_n})$ converge to $p(F,Q_{\V x})$.
\end{remark}

\begin{lemma}\label{lm:compact}
The space $\PC$ and distance $\dedit$ have the following properties.
\begin{enumerate} 
  \item[(i)]\label{it:compact} The space $(\PC,\dedit)$ is a compact metric space.
\item[(ii)]\label{it:dense} The set 
of complete partite graphs $\CP$ is dense in $(\PC,\dedit)$.
 \item[(iii)]\label{it:pCts} The function $\lambda$  can be extended to a
continuous function on the whole of $\PC$, namely by defining
\begin{align*}
\lambda(\V x) &:= \lim_{n \to \infty}\lambda(G_{n,\V x}),\quad\text{for }x \in \PC.
\end{align*}
\end{enumerate}\end{lemma}

\bpf
We begin with (i). From the definitions it is clear that $\dedit(\V x, \V y)=\dedit(\V y, \V x)$ for all $\V x, \V y \in \PC$.
By Lemma~\ref{lm:edit}(iii), $\dedit$ defined on $\PC$ satisfies the triangle inequality.
Finally, by definition, $\dedit(\V x,\V y)=0$ if and only if $\V x = \V y$.
So $(\PC, \dedit)$ is a metric space.
To show that it is compact, Proposition~\ref{lem-edit-pw} implies that it suffices to show that $\PC$ is compact under the topology of pointwise convergence.
For this, let $(\V x_n)_n$ be an infinite sequence of elements of $\PC$. Then we can define
its accumulation point $\V y$ iteratively as follows. Initially let $i=0$. By passing to
a subsequence of $(\V x_n)_n$, we may assume that $(x_{n,i+1})_n$ converges to some $y_{i+1} \in \mathbb{R}$. If 
$y_{i+1}=0$, then stop and output $\V y := (y_1,\ldots,y_i, 
0,0,\ldots)$. Otherwise,
increase $i$ by one and continue. If the iteration does not terminate, output $\V y := (y_1,y_2,\ldots)$.
One can easily see that $\V y$ is indeed an accumulation point of $(\V x_n)_n$, completing the proof of (i).
Alternatively, the compactness of $\PC$ follows from observing that $\PC$ is a closed subset of the compact space $[0,1]^{\mathbb{N}}$.

Part~(ii) immediately follows since for every $\V x\in\PC$, the sequence $(G_{n,\V x})_n$ of complete partite graphs converges in edit distance to $\V x$.
Indeed, for each $n \in \mathbb{N}$ we have that $\V x_{G_{n, \V x}} \in \CP$, and the definitions imply that $\V x_{G_{n,\V x}}$ converges pointwise to $\V x$.
By Proposition~\ref{lem-edit-pw}, it also converges in edit distance.

It remains to prove~(iii).
Recall that we fixed a function $\gamma : \mathcal{G}_k \rightarrow \mathbb{R}$ and for all $n \in \mathbb{N}$ and $G \in \mathcal{G}_n$, we defined $\lambda(G)$ as in~(\ref{eq:gamma}).
Let $\gammamax := \max\{|\gamma(F)|: F \in \mathcal{G}_k\}$ (which exists since the domain of $\gamma$ is finite for fixed $k$).
Let $n \in \mathbb{N}$, let $G \in \mathcal{G}_n$ and let $xy$ be a pair in $V(G)$. Then
$$
|\lambda(G) - \lambda(G \oplus xy)| \leq \binom{n}{k}^{-1} \sum_{\stackrel{X \in \binom{V(G)}{k}:}{x,y \in X}} [\gamma(G[X])-\gamma((G\oplus xy)[X])]\leq \frac{\binom{n-2}{k-2}\cdot 2\gamma_{\mathrm{max}}}{\binom{n}{k}} = \frac{2\binom{k}{2}\gamma_{\mathrm{max}}}{\binom{n}{2}}.
$$
Therefore, using the triangle inequality, we have for any $G,H \in \mathcal{G}_n$ that
\begin{equation}\label{lipschitz}
|\lambda(G)-\lambda(H)| \leq 2\binom{k}{2}\gamma_{\mathrm{max}}\cdot \done(G,H) + O(1/n) \leq 6\binom{k}{2}\gamma_{\mathrm{max}}\cdot \dedit(G,H) +O(1/n),
\end{equation}
where the final inequality follows from Lemma~\ref{lm:edit}(i).
Thus 
$$
|\lambda(G_{n,\V x})-\lambda(G_{n,\V y})| \leq 6\binom{k}{2}\gamma_{\mathrm{max}} \cdot \dedit(G_{n, \V x}, G_{n,\V y})+O(1/n),
$$
and by (i) we have that the function $\lambda : \PC \rightarrow \I R$ given by
$
\lambda(\V x) := \lim_{n\rightarrow \infty}\lambda(G_{n,\V x})
$
is well-defined for all $\V x \in \PC$ and is continuous with respect to~$\dedit$.
\epf

Note that the extension in Part~(iii) of Lemma~\ref{lm:compact} is unique since  $\CP$ is dense in $\PC$. 

The lemma implies that 
$\lambdamax:=\lim_{n\to\infty}\lambda(n)$ defined in the introduction can equivalently be defined as $\lambdamax := \max\{\lambda(\V x): \V x \in \PC\}$.
Moreover,
for every $\V x=(x_1,x_2,\ldots)\in\PC$, we have that $\lambda(\V x)$ has the 
 following analytic formula.
 Let $\omega_1,\ldots,\omega_k$ be independent samples from
 $\Omega_{\V x}$ which is the probability space on 
  $\I N_0:=\{0,1,2,\ldots\}$
 where the probability of $i$ is $x_i$. Let the \emph{random sample} $\I G(\V 
 x,k)$ be
 equal to 
  $$
   G_{\V x}^{\omega_1,\ldots,\omega_k}
  :=\left([k],{[k]\choose 2}\setminus \left(\bigcup_{i\in\I N} {\{j\mid 
  \omega_j=i\}\choose 
  2}\right)\right),$$
   which is the complete graph 
  on $[k]$ except we do not connect two distinct indices $j,h\in [k]$
 if $\omega_j=\omega_h\not=0$. 
 One can show using the Chernoff bound and the Borel-Cantelli lemma that
 $(\I G(\V x,n))_{n}$ converges to $\V x$ in $\PC$ with probability 1
(see e.g.~the more general Proposition~11.32 in~\cite{Lovasz12}).
Clearly we have that
$$
\lambda(\V x) = \I E (\gamma(\I G(\V x, k))).
$$

We let $\OPT$  consist of all maximisers $\V 
x\in\PC$, that is,
\begin{equation*}
\OPT = \OPT(\lambda) := \lbrace \V x \in \PC : \lambda(\V x) \geq \lambda(\V y) \text{ for all } \V y \in \PC\rbrace = \{ \V x \in \PC: \lambda(\V x)=\lambdamax \} \neq \varnothing.
\end{equation*}
The non-emptiness assertion follows from Lemma~\ref{lm:compact}(i) and~(iii).
Let us see why the forward inclusion of the third equality is true. Take any $\V x\in\PC$ with $\lambda(\V x) \geq \lambda(\V y)$ for all $\V y \in \PC$.
For each $n \in \mathbb{N}$, since $\lambda$ is symmetrisable, there is a complete partite graph $F_n$ on $n$ vertices such that $\lambda(F_n)=\lambda(n)$. 
Let $\V y_n := \V x_{F_n}$.
For any $\V x \in \PC$, we have $\lambda(G_{n, \V x}) \leq \lambda(F_n) = \lambda(n)$.
By passing to a subsequence we may assume that $\V y_{n_i}$ converges to some~$\V y\in\PC$.
Then $\lambda(\V x) \leq \lambda(\V y) = \lim_{i\to\infty}\lambda(n_i)=\lambdamax$. Thus we must have $\lambda(\V x) = \lambdamax$, as desired.

This definition of $\OPT$ is equivalent to the one in the introduction. Indeed, 
let $\V a = (a_1,a_2,\ldots) \in \PC$
be such that there exists a sequence $(H_n)_n$ of complete partite graphs such that, as $n\to\infty$, we have $v(H_n) \rightarrow \infty$, $\lambda(H_n)\to \lambdamax$ and for every $i\geq 1$ the number of vertices in the $i$-th largest part of $H_n$ is $(a_i+o(1))v(H_n)$.
Then $\V x_{H_n}\to\V a$ and $\lambda(\V a) = \lim_{n\to\infty}\lambda(H_n) = \lambdamax$, as required.
On the other hand, let $\V x \in \PC$ be such that $\lambda(\V x)=\lambdamax$.
Then $(G_{n,\V x})_n$ is the required sequence of graphs.

\subsection{Polynomials}

We will be interested in various functions on $\PC$, in particular the extension of $\lambda$ from the family of complete partite graphs to $\PC$.
For these we introduce a notion of polynomial on $\PC$ which will help us prove that functions related to $\lambda$ are continuous.

Let $\Sigma(d) := \{(d_1,\ldots,d_t) \in \mathbb{N}^t: t \in \mathbb{N}_0\text{ and }d_1+\ldots+d_t=d\}$ be the set of ordered tuples of positive integers summing to $d$.
Let $S_{\varnothing}(\V x) := 1$ and for $t \in \mathbb{N}$ and $\V d := (d_1,\ldots,d_t) \in \Sigma(d)$, define an \emph{elementary symmetric polynomial} $S_{\V d} : \{\V x \in \mathbb{R}^{\mathbb{N}}: \|\V x\|_1 < \infty\} \to \mathbb{R}$ by
\begin{equation}\label{eq:S}
S_{\V d}(\V x) = S_{d_1,\ldots,d_t}(\V x):=\sum_{\substack{\mathrm{distinct}\\i_1,\ldots,i_t\in \I N}} \prod_{j=1}^t x_{i_j}^{d_j}.
\end{equation}
Since $\sum_{\substack{\mathrm{distinct}\\i_1,\ldots,i_t\in \I N}} \prod_{j=1}^t |x_{i_j}|^{d_j} \leq \left(\sum_{i \geq 1}|x_i|\right)^{d} <\infty$,
each $S_{\V d}(\V x)$ converges absolutely.

We say that a function $p: \PC \rightarrow \mathbb{R}$ is a \emph{$\PC$-polynomial} 
if it can be written as a finite polynomial of $S_{\V d}^I(\V x) := S_{\V d}(\V x^I)$ for $I \subseteq \mathbb{N}$, where $\V x^I \in \PC$ is obtained from $\V x$ by removing every $x_i$ with $i \in I$ and moving back remaining entries to fill in the `gaps'. (Thus  $S_{\V d}^I(\V x)$ is defined by the version of~\eqref{eq:S} where the sum is restricted to indices not in~$I$.)
So, for example, $x_0 = S_{\varnothing}(\V x) - S_{1}(\V x)$, $x_i=S_1(\V x) - S_1^{\{i\}}(\V x)$ for $i \in \mathbb{N}$ and $x_1+x_3+x_5+x_7+\ldots=S_1^{2\cdot\mathbb{N}}(\V x)$ are $\PC$-polynomials, while $x_1+2x_2+3x_3+\ldots$ is not.
Given any $\PC$-polynomial $p$, there is a finite partition $\mathbb{N}=I_1 \cup \ldots \cup I_s$
such that $p(x_1,x_2,\dots) = p(y_1,y_2,\ldots)$ where $\V y$ is any element of $\PC$ obtained from $\V x$ by permuting indices within each part $I_i$.
Indeed, one can obtain $I_1,\ldots,I_s$ by grouping together indices that belong to exactly the same sets $I$ in the definition of $p$.

Take any $m\in\mathbb N$ and $\V d = (d_1,\ldots,d_t) \in \mathbb{N}^t$. Consider $\frac{1}{h}(S_{\V d}(\V x')-S_{\V d}(\V x))$ where, for all $i \geq 1$ we have $x_i'=x_i$, except $x_m'=x_m+h$, and let $h \to 0$.
Apply the binomial expansion to each $(x_m+h)^{d_j}$. As all series converge absolutely, we can change the order of summation and collect the same powers of $h$. We obtain
$$
\frac{S_{\V d}(\V x')-S_{\V d}(\V x)}{h} = \sum_{\substack{\mathrm{distinct}\\i_1,\ldots,i_t\in \I N}} \frac{\partial}{\partial x_m}\left(\prod_{j=1}^t x_{i_j}^{d_j}\right) + \delta
$$
where $\delta$ is an error term satisfying $|\delta| \leq h\cdot 2^d$.
So
 we can define partial derivatives $\frac{\partial p}{\partial x_i}$ for $i=1,2,\ldots$ via term-by-term differentiation.
Also, if $p=s(S_{\V d} : \V d \in \mathbb{N}^{\leq k})$ where $s$ is a finite polynomial, then define $\frac{\partial p}{\partial x_0} := -\frac{\partial s}{\partial S_{1}}$.
Thus we can define partial derivatives of any $\PC$-polynomial, and each such derivative is itself a $\PC$-polynomial.
For a complete partite graph $G$ on $n$ vertices with parts $V_1,\ldots,V_m$ of size at least $2$ and clique part $V_0$, define for $I \subseteq \mathbb{N}$
\begin{align}
S_{\V d}^I(G)&:=(d!)^{-1}\binom{n}{d}^{-1}\sum_{\substack{\text{distinct }\\ i_1,\ldots,i_t \in [m]\setminus I}}\prod_{1 \leq j \leq t}d_j!\binom{|V_{i_j}|}{d_j}
\label{Sdgraph}&= \sum_{\substack{\text{distinct }\\i_1,\ldots,i_t \in [m]\setminus I}}\prod_{1 \leq j \leq t}\left(\frac{|V_{i_j}|+O(1)}{n}\right)^{d_j},
\end{align}
and let $S_{\V d}(G) := S_{\V d}^\varnothing(G)$.
So $S_{\V d}^I(G)$ is equal to $S_{\V d}(G^I)$, up to a scaling factor, where $G^I:= G-\bigcup_{i \in [m]\cap I}V_i$.

\begin{lemma}\label{poly}
Let $d$ be an integer and let $\V d = (d_1,\ldots,d_t) \in \Sigma(d)$. Then
\begin{itemize}
\item[(i)] $S_{\V d}$ is uniformly continuous on $(\PC,\dedit)$.
\item[(ii)] Each $\PC$-polynomial is uniformly continuous on $(\PC,\dedit)$.
\item[(iii)] For all $\V x \in \PC$ we have $S_{\V d}(\V x) = \lim_{n\to \infty}S_{\V d}(G_{n,\V x})$.
\end{itemize}
\end{lemma}

\bpf
We start with (i). By Proposition~\ref{lem-edit-pw}, convergence in edit distance and pointwise convergence induce the same topology on $\PC$. By Lemma~\ref{lm:compact}(i), $\PC$ is compact.
Therefore it suffices to show that each $S_{\V d}$ is continuous under pointwise convergence, which is e.g.~given by the metric $d(\V x, \V y) := \sum_{i \geq 1}2^{-i}|x_i-y_i|$.
For this, let $\eps>0$ and let $\delta = 2^{-8d/\eps}$.
Let $\V x, \V y \in \PC$ satisfy $d(\V x, \V y) \leq \delta$.
Choose $M  = \lceil\log_2(\delta^{-1/2})\rceil$ (so $1/M \leq \eps/(4d)$) and
let $\V x' = (x_1,\ldots,x_{M},0,\ldots)$ and $\V y' = (y_1,\ldots,y_M,0,\ldots)$.
Then $d(\V x, \V x') = \sum_{i > M}2^{-i}x_i \leq 2^{-M} \leq \sqrt{\delta}$.
So $d(\V x', \V y') \leq 3\sqrt{\delta}$.
Moreover,
$$
S_{\V d}(\V x)-S_{\V d}(\V x') = 
\sum_{1 \leq s \leq t}\sum_{i>M}x_i^{d_s}S_{\V d^{(s)}}(\V x^{(i)})
\leq tx_{M+1} \leq d/M \leq \eps/4,
$$
where $\V d^{(s)}= (d_1,\ldots,d_{s-1},d_{s+1},\ldots,d_t)$
and $\V x^{(i)} = (x_1,\ldots,x_{i-1},x_{i+1},\ldots)$.
Similarly $S_{\V d}(\V y)-S_{\V d}(\V y') \leq \eps/4$.
Now, $S_{\V d}(\V x')$ is a polynomial in at most $M$ variables.
For each $1 \leq i \leq M+1$, let $\V z_i := (x_1,\ldots,x_{i-1},y_i,y_{i+1},\ldots,y_M,0,\ldots)$. Then
$$
|S_{\V d}(\V x') - S_{\V d}(\V y')| = |S_{\V d}(\V z_1)-S_{\V d}(\V z_{M+1})| \leq \sum_{i=1}^{M} |S_{\V d}(\V z_{i+1})-S_{\V d}(\V z_{i})|.
$$
Now
$$
S_{\V d}(\V z_{i+1})-S_{\V d}(\V z_{i}) = \sum_{1 \leq s \leq t}(x_i^{d_s}-y_i^{d_s})S_{\V d^{(s)}}(\V z_i^{(i)})
=p_i(x_i)-p_i(y_i)
$$
where we view $p_i$ as a polynomial in one variable.
Thus $p_i$ is Lipschitz with constant at most $\max_{z \in [0,1]}|p_i'(z)| \leq d_1+\ldots+d_t = d$.
So $|p_i(x_i)-p_i(y_i)| \leq d|x_i-y_i|$.
Thus
$$
|S_{\V d}(\V x)-S_{\V d}(\V y)| \leq \eps/2 + d\sum_{i=1}^{M}|x_i-y_i| \leq \eps/2 + d2^M d(\V x, \V y) \leq \eps/2+\sqrt{\delta} < 2\eps/3,
$$
completing the proof of~(i).

Now (ii) follows immediately since every $S_{\V d}$ is bounded, and sums and products of bounded uniformly continuous functions are uniformly continuous.

For~(iii), fix $\V x \in \PC$.
In $G_{n,\V x}$, writing $V_i^n$ for the $i$th part, we have each $(|V_i^n|+O(1))/n \to x_i$ as $n\to\infty$, so as $S_{\V d}$ is continuous we have $S_{\V d}(G_{n,\V x}) \to S_{\V d}(\V x)$.
\epf

\section{Strictness and a restatement of the main result}\label{sec-strict}

In this section, we will finally define what it means for $\lambda$ to be `strict'.
Very roughly speaking, it means that when an elementary change is made to a complete partite graph on which $\lambda$ is maximised, the decrease in $\lambda$ is as much as it possibly could be.
An `elementary change' is either `flipping a pair' (changing a non-edge to an edge or vice versa); or adding a vertex which is either adjacent to every vertex in a part, or to no vertex in a part.
It seems that it is more convenient to state this property in terms of limits rather than graphs (which is why the definition is deferred until now). 
We will first make the relevant definition and then discuss it further.

\subsection{Definitions and notation}

\begin{definition}[$\flip_{i_1i_2}\lambda$ and $\join_{b,\alpha}\lambda$]\label{def-nabla}
Given an $n$-vertex graph $G=(V,E)$ and a pair $x,y$ of vertices of $G$, define
$$
\flip_{xy}\lambda(G) := \frac{1}{{n-2\choose k-2}} (\Lambda(G)-\Lambda(G\oplus xy)).
$$
Given $\V x \in \PC$ and $i_1,i_2 \in \supp^*(\V x)$, define
 $$
 \flip_{i_1i_2}\lambda(\V x) := \lim_{n \rightarrow \infty}\flip_{v_1v_2}\lambda(G_{n,\V x}),
 $$
 where $v_1$, $v_2$ are distinct vertices of the vertex classes $V_{i_1}$ and $V_{i_2}$ of $G_{n,\V x}$ respectively.

For all $i \in \mathbb{N}_0$, we define $e_i$ to be the function $e_i : \mathbb{N}\rightarrow \lbrace 0,1\rbrace$ with $e_i(j)=0$ if and only if $j=i$ (so $e_0 \equiv 1$).
Let $b : \mathbb{N}\rightarrow \lbrace 0,1\rbrace$ and $\alpha \in [0,1]$. 
We write $G+_{b,\alpha} u$ for the graph obtained from $G$ with $\PC$-structure $V_0,V_1,\ldots,V_m$ by adding a new vertex $u$ and, for $i \geq 1$, adding every edge between $u$ and $V_i$ if $b(i)=1$, and no edges otherwise; and adding $\lfloor \alpha |V_0|\rfloor$ edges between $u$ and $V_0$.
Define
$$
\join_{b,\alpha}\lambda(G) := \frac{1}{\binom{n}{k-1}}
(\Lambda(G +_{e_1,1} u,u) - \Lambda(G +_{b,\alpha} u,u))
=\frac{1}{\binom{n}{k-1}}
(\Lambda(G +_{e_1,1} u) - \Lambda(G +_{b,\alpha} u))
$$ 
where $u \notin V(G)$, and let 
\begin{align*}
\join_{b,\alpha}\lambda(\V x) := \lim_{n\to\infty}\join_{b,\alpha}\lambda(G_{n,\V x})\quad
\text{and}\quad
\lambda(\V x, (b,\alpha)) := \lim_{n\rightarrow\infty}\lambda(G_{n,\V x}+_{b,\alpha} u,u).
\end{align*}
By convention take $\alpha=1$ if $x_0=0$ (when $V_0=\varnothing$). 
\end{definition}

Given $k_0 \in \mathbb{N}_0$ and a tuple $\V k = (k_1,\ldots,k_t)$ of positive integers, define the graph $G_{\V k}^{k_0}$ as follows.
Let $G_{\V k}^{k_0}$ be the complete partite graph with $t$ parts $U_1,\ldots,U_t$ of size $k_{1},\ldots,k_{t}$ respectively, together with an additional $k_0$ singletons $x_1,\ldots,x_{k_0}$, whose union is denoted by $U_0$.

Both limits in Definition~\ref{def-nabla} exist and each $\lambda,\flip_{i_1i_2}\lambda$, $\join_{b,\alpha}\lambda$, $\lambda(\cdot,(b,\alpha))$ is a $\PC$-polynomial.
Indeed, since each $G_{n,\V x}$ is a complete partite graph (with parts $V_0^n,V_1^n,\ldots$), the quantities $\lambda(G_{n,\V x}), \flip_{v_1v_2}\lambda(G_{n,\V x})$ and $\join_{b,\alpha}\lambda(G_{n,\V x})$ are finite polynomials in variables $|V_0^n|$ and $S_{\V d}^I(G_{n,\V x})$ for $\V d \in \Sigma(d)$ with $d \leq k$ and $I \subseteq \mathbb{N}$.
Indeed, for $\lambda$ we need only $I=\varnothing$;
for $\flip_{v_1v_2}\lambda$ we could take only $I=\varnothing,\{v_1\},\{v_2\},\{v_1,v_2\}$ and their complements,
and for $\join_{b,\alpha}\lambda$, $I=\varnothing,\ \supp(b)$ and their complements.
Thus, by Lemma~\ref{poly}, $\flip_{i_1i_2}\lambda$ and $\join_{b,\alpha}\lambda$ are $\PC$-polynomials.

In fact, one can explicitly write these polynomials.
For positive integers $b_1 \geq \ldots \geq b_r$, let ${\rm sym}(b_1,\ldots,b_r)$ be the number of permutations of $[r]$ that keep the sequence $(b_1,\ldots,b_r)$ unchanged. In other words, if we take $i_0 := 1 < i_1 < \ldots < i_q < r+1 =: i_{q+1}$ such that $b_i=b_{i'}$ if and only if there is $j \in [q]$ such that $i_{j-1} \leq i,i' < i_{j}$ then ${\rm sym}(b_1,\ldots,b_r)= (i_1-i_0)!\,\ldots\, (i_{q+1}-i_{q})!$.
Also, write $\binom{t}{t_1,\ldots,t_s}:=t!(t_1!\,\ldots\, t_s!)^{-1}$ when $\sum_{i=1}^{s}t_i=t$.
Consider $p(K_{a_1,\ldots,a_\ell},\cdot)$, which is one instance of $\lambda$, where $a_1,\ldots,a_\ell$ are in non-increasing order, and let $t \in [\ell]$ be the largest integer such that $a_1,\ldots,a_t \geq 2$. 
Then we have the following analytic formula:
\begin{equation}\label{eq:lambdaan}
p(K_{a_1,\ldots,a_\ell},\V x) = p(K_{a_1,\ldots,a_\ell},Q_{\V x}) = \frac{\binom{a_1+\dots+a_\ell}{a_1,\ldots,a_\ell}}{{\rm sym}(a_1,\ldots,a_{\ell})}\sum_{0 \leq s \leq \ell-t}\binom{\ell-t}{s}x_0^s\cdot S_{a_1,\ldots,a_{\ell-s}}(\V x).
\end{equation}
Using~(\ref{eq:lambdaan}), one can write $\flip_{i_1i_2}\lambda$ and $\join_{b,\alpha}\lambda$ as explicit $\PC$-polynomials.


The next proposition gives that for all $\V x \in \OPT$, $\join_{e_i,1}\lambda(\V x) = 0$ for all $i \in \supp^*(\V x)$, which corresponds to saying that every vertex in the realisation of an optimal $\V x$ contributes optimally to $\lambda$. Thus
$$
  \join_{b,\alpha}\lambda(\V x)=\lambda(\V x)-\lambda(\V x, (b,\alpha)) \quad\text{for all }(b,\alpha) \text{ and }\V x \in \OPT.
$$

\begin{proposition}\label{lagmult}
Define $k$ and $\lambda$ as in~(\ref{eq:gamma}). The following hold for all $\V x \in \PC$.
\begin{itemize}
\item[(i)] For all $i \in \supp^*(\V x)$ we have
$
\frac{1}{k}\cdot\frac{\partial\lambda(\V x)}{\partial x_i} = \lambda(\V x,(e_i,1))$.
\item[(ii)] If in addition $\V x \in \OPT$, then for all $i \in \supp^*(\V x)$ we have
$\frac{1}{k}\cdot \frac{\partial\lambda(\V x)}{\partial x_i}=\lambda(\V x)$.
\item[(iii)] The following pairs differ by $O(1/n)$ as $n \to \infty$: $\{\lambdamax, \lambda(n)\}$, $\{\lambda(\V x), \lambda(G_{n,\V x})\}$ and $\{\lambda(G_{n,\V x},u), \lambdamax\}$, the last pair for $\V x \in \OPT$ and $u \in V(G_{n,\V x})$.
\end{itemize}
\end{proposition}

\bpf
The equality in (i)~can be checked directly.

For (ii),
the theory of Lagrange multipliers implies that, for all $i,j \in \supp^*(\V x)$, we have $\frac{\partial\lambda(\V x)}{\partial x_i} = \frac{\partial\lambda(\V x)}{\partial x_j}$.
Indeed, if we fix the rest of $\V x$ apart from $x_i,x_j$, fix $s=x_i+x_j$ and vary $x_i,x_j$, then we can view $\lambda$ as a polynomial in $x_i,x_j$ (of degree at most $k$). Introducing a new variable $\mu$, the Lagrangian is $\mathcal{L}(x_i,x_j,\mu) = \lambda(\V x)-\mu(x_i+x_j-s)$. 
The stationary points of $\mathcal{L}$ occur when $(\frac{\partial\mathcal{L}}{\partial x_i},\frac{\partial\mathcal{L}}{\partial x_j},\frac{\partial\mathcal{L}}{\partial \mu})=(0,0,0)$, i.e.\ when $\frac{\partial\lambda(\V x)}{\partial x_i}-\mu=\frac{\partial\lambda(\V x)}{\partial x_j}-\mu$, as required.
Since $\lambda$ is a $\PC$-polynomial with each monomial having total degree $k$,
we have for all $i \in \supp^*(\V x)$ that
$$
\frac{\partial\lambda(\V x)}{\partial x_i} =\sum_{j \geq 0}x_j\frac{\partial\lambda(\V x)}{\partial x_j}  = k\cdot\lambda(\V x),
$$
giving the required.

Let us turn to Part~(iii). The inequality $|\lambdamax-\lambda(n)|=O(1/n)$ follows from a standard blow-up trick, see e.g.~\cite[Lemma~2.2]{PikhurkoSliacanTyros19}. The claim for the second pair
follows from the fact that each named function on $\PC$ is a $\PC$-polynomial, a finite polynomial of $S_{\V d}(G^I_{n,\V x})$ terms, so the error bound comes from~(\ref{Sdgraph}) when applied to $G_{n, \V x}$.
For the last claim of Part~(iii), a version of~(\ref{Sdgraph}) implies that $|\lambda(G_{n,\V x},u)-\frac{1}{k}\cdot \frac{\partial\lambda(\V x)}{\partial x_i}|=O(1/n)$ where $u$ is in the $i$-th part of $G_{n,\V x}$. Then Part~(ii) gives the required.
\epf

\begin{corollary}\label{unicont}
For every $\eps>0$ there exists $\delta>0$ such for all $\V x,\V y \in \PC$ with 
$\dedit(\V x, \V y) \leq \delta$, we have
$$
|\lambda(\V x)-\lambda(\V y)|, |\join_{b,\alpha}(\V x)-\join_{b,\alpha}(\V y)|, |\lambda(\V x,(b,\alpha))-\lambda(\V y,(b,\alpha))| \leq \eps
$$
for all $b : \mathbb{N}\rightarrow \lbrace 0,1\rbrace$ and $0 \leq \alpha \leq 1$, and
$$
|\flip_{i_1i_2}\lambda(\V x)-\flip_{i_1i_2}\lambda(\V y)| \leq \eps
$$
for all $i_1,i_2 \in \supp^*(\V x) \cap \supp^*(\V y)$.
\end{corollary}

\bpf
We have seen that each function $\lambda,\flip_{i_1i_2}\lambda,\join_{b,\alpha}\lambda,\lambda(\cdot,(b,\alpha))$ is a $\PC$-polynomial with degree at most $k$ and with coefficients whose absolute values are bounded. Thus Lemma~\ref{poly} implies that the family of $\lambda,\flip_{i_1i_2}\lambda,\join_{b,\alpha}\lambda,\lambda(\cdot,(b,\alpha))$ over all $i_1,i_2,b,\alpha$ is uniformly equicontinuous, as required.
\epf

The following crucial definition of the \emph{strictness} property of a function $\lambda$ requires that both $\flip_{i_1i_2}\lambda$ and $\join_{b,\alpha}\lambda$ are bounded from below whenever $(b,\alpha)$ is not close to some $(e_i,1)$. Roughly speaking, this means that $\lambda$ is sensitive to small alterations in a graph.

\begin{definition}[Strictness]\label{defn:strict}
	We say that $\lambda$ is \emph{strict (with parameter $c$)} if there is $c=c(\lambda)>0$ such that for each $\V x\in \OPT$, we have
	\begin{enumerate} 
	\item[\emph{\Sone}]\label{itm-S1} $\flip_{i_1i_2}\lambda(\V x) \geq c$ for all $i_1,i_2 \in \supp^*(\V x)$,
	\item[\emph{\Stwo}]\label{itm-S2} $\join_{b,\alpha}\lambda(\V x) \geq c((1-\alpha)x_0+ \min_{i \in \supp^*(\V x)}w_i)$, where
$$
w_i := \mathds{1}_{i > 0}b_ix_i + \sum_{j\in\supp^*(\V x)\setminus \lbrace 0,i\rbrace}(1-b_j)x_j.
$$
	\end{enumerate}
\end{definition}


In the next two subsections, we will motivate these definitions, which appear somewhat complicated at first sight.

\subsubsection{$\flip_{i_1i_2}\lambda$: flipping a pair of vertices}
Take a complete partite graph $G$ of large order $n$ such that $\lambda(G)\approx \lambdamax$ and let $G'=G\oplus xy$ be obtained by flipping the adjacency of an arbitrary pair $xy\in{V\choose 2}$.
Then the number of vertex subsets of size $k$ which contain both $x$ and $y$ is $\binom{n-2}{k-2}$, so in the worst case, $\gamma$ decreases by a constant for all such subsets, and thus $\lambda$ decreases by $\Omega(\binom{n-2}{k-2}/\binom{n}{k}) = \Omega(1/n^2)$.
Property~\Sone~says that this worst-case behaviour is realised for every `wrong' pair $xy$.

Observe that
\begin{equation*}
 \flip_{i_1i_2}\lambda(\V x)= \mathbb{E}_{\omega_3,\ldots,\omega_k\sim \Omega_{\V x}} 
 [\gamma(G_{\V 
 x}^{i_1,i_2,\omega_3,\ldots,\omega_k})-\gamma(G_{\V x}^{i_1,i_2,\omega_3,\ldots,\omega_k}\oplus \{1,2\})],
 \end{equation*}
  that is, we look at the conditional expectation of the change 
 in $\lambda$ if we flip the pair $\{1,2\}$ in
a random sample $\I G(\V x,k)$ conditioned on $\omega_1=i_1$ and 
$\omega_2=i_2$.

\subsubsection{$\join_{b,\alpha}\lambda$: adding a new vertex}

Again consider a complete partite graph $G$ of large order $n$ such that $\lambda(G) \approx \lambdamax$ and obtain a graph $G'$ from $G$ by adding a new vertex $u$ which, for each part of $G$, either connects to all or none of its vertices (here we are thinking of $V_0$, if it exists, as consisting of $|V_0|$ singleton parts).
If the attachment of $u$ mirrors an existing vertex, then its contribution to $\lambda$ is approximately $\lambdamax$ (and $G'$ is the same as $G$ in the limit).
But, if not, as $u$ lies in $\binom{n}{k-1}$ subsets, in the worst case, $\lambda$ decreases by $\Omega(\binom{n}{k-1}/\binom{n+1}{k})=\Omega(1/n)$.
Property~\Stwo~says that this worst-case behaviour is realised for every $u$ with `wrong' attachment.

Suppose that $G_{n,\V x}$ has $\PC$-structure $V_0,V_1,\ldots,V_{m(n)}$. Then, for $0 \leq i \leq m(n)$, let
$W_i$ be the minimum number of edits needed to move the vertex $u$ in $G_{n,\V x} +_{b,\alpha} u$ into the $i$-th part.
So each $W_i$ being large corresponds to $u$ being attached in an atypical manner, and some $W_i$ small means that $u$ behaves like an existing vertex.
It is not hard to show that
$
\lim_{n\rightarrow\infty}W_i/n = w_i+(1-\alpha)x_0
$, and of course if $b=e_i$ and $\alpha=1$, then $w_i+(1-\alpha)x_0=0$ (since no edits are needed to move $u$ to the $i$-th part).
So~\Stwo~requires that, whenever $n$ is large, the contribution to $\lambda$ lost by a vertex $u$ in $G_{n,\V x} +_{b,\alpha} u$ is a significant fraction of the number of edits needed to fit $u$ into $G_{n,\V x}$.

Observe that (using Proposition~\ref{lagmult} and the remark immediately before it)
\begin{equation*}
 \join_{b,\alpha}\lambda(\V x):= \mathbb{E}_{\omega_1,\ldots,\omega_k\sim \Omega_{\V x}} 
 [\gamma(G_{\V 
 	x}^{\omega_1,\ldots,\omega_k})]-
 \mathbb{E}_{\omega_1,\ldots,\omega_{k-1}\sim \Omega_{\V x}} [\gamma(G_{\V x}^{\omega_1,\ldots,\omega_{k-1}}+_{b,\alpha}u)],
\end{equation*}
where $G_{\V x}^{\omega_1,\ldots,\omega_{k-1}}+_{b,\alpha}u$ is the random graph obtained by adding $u$ to $\I G(\V x, k-1)$ with $ui$ an edge when $\omega_i \neq 0$ if and only if $b(\omega_i)=1$, and $ui$ an edge when $\omega_i=0$ with probability $\alpha$.

\subsection{Main result}

We are now ready to precisely state the `limit version' of our main result.

\begin{theorem}\label{thm-main}
Let $k$ be a positive integer and let $\gamma : \mathcal{G}_k \rightarrow \mathbb{R}$. Define $\lambda : \mathcal{G} \rightarrow \mathbb{R}$ by setting $\lambda(G):= \binom{n}{k}^{-1}\sum_{X \in \binom{V}{k}}\gamma(G[X])$ for all $G \in \mathcal{G}_n$ and $n \in \mathbb{N}$, and let $\lambda(n) := \max_{G \in \mathcal{G}_n}\lambda(G)$.
Suppose that $\lambda$ is symmetrisable and $|\OPT(\lambda)|<\infty$.
	Then $\lambda$ has perfect stability if it is strict.
\end{theorem}

The following corollary states that strict symmetrisable functions exhibit classical stability, in the sense that any sufficiently large graph which is sufficiently close to being optimal can be edited by changing an arbitrarily small fraction of its adjacencies to obtain a complete partite graph with the correct part sizes.

\begin{corollary}\label{cor-main}
Define $k$ and $\lambda$ as in~(\ref{eq:gamma}) and suppose that they satisfy the assumptions in Theorem~\ref{thm-main}, and suppose further that $\lambda$ is strict.
Then for all $\eps>0$ there exist $\delta,n_0>0$ such that for every graph $G$ of order $n \geq n_0$ for which $\lambda(G) \geq \lambdamax-\delta$, there is $\V x \in \OPT(\lambda)$ for which $\dedit(G,\V x) \leq \eps$.
\end{corollary}

\bpf
Let $c = c(\lambda)>0$ be such that $\lambda$ is strict with parameter $c$.
Apply Theorem~\ref{thm-main} to obtain $C$ such that $\lambda$ is perfectly stable with constant $C$.
Suppose that the statement does not hold. Then there is a sequence of counterexamples $(G_n)_{n}$ with $v_n := v(G_n)\rightarrow\infty$ such that
$\lambda(G_n) \geq \lambdamax - 1/n$ but $\dedit(G_n,\V x) > \eps$ for all $\V x \in \OPT$.
By taking a subsequence if necessary, we may assume that each $v_n \geq n$. Let $n$ be sufficiently large.
By Theorem~\ref{thm-main}, there is some $H_n \in \CP_{v_n}$ for which
$$
\done(G_n,H_n)/C \leq \lambda(v_n)-\lambda(G_n) \leq \lambda(v_n)-\lambdamax +O(1/v_n)  \leq O(1/n)
$$ where we used Proposition~\ref{lagmult}(iii). But then by~(\ref{lipschitz}),
$$
\lambda(H_n) \geq \lambda(G_n) - 2\binom{k}{2}\gammamax \done(G_n,H_n)-O(1/v_n) \geq \lambdamax - O(1/n)\left(1 + 2C\binom{k}{2}\right).
$$
So, writing $\V x_n := \V x_{H_n}$, and taking a subsequence if necessary, we see that $\V x_n \rightarrow \V x \in \OPT$.
But then, when $n$ is sufficiently large, using Lemma~\ref{lm:edit},
$$
\dedit(G_n,\V x) \leq \dedit(G_n,H_n) + \dedit(\V x_n, \V x) \leq \done(G_n,H_n)+\dedit(\V x_n, \V x) < \eps,
$$
a contradiction.
\epf

\section{Finitely many maximisers}\label{se:ff}

We will need the following result which states that if the limit problem has finitely many optimisers, then all non-zero entries in them are separated from $0$ by some constant $\beta>0$.

\begin{lemma}\label{lm:ff} If $|\OPT|<\infty$ then there is $\beta>0$ such that $\OPT \subseteq \PC_\beta$.
\end{lemma}

The rest of the section is dedicated to proving Lemma~\ref{lm:ff}. Our proof is an adaptation of the proof of Glebov, Grzesik, Klimo\v sov\'a and Kr\'al'~\cite{GlebovGrzesikKlimosovaKral15jctb} who, in particular, worked on the finite forcibility of graphons which are a countable union of cliques.
Recall notions related to graphons in Section~\ref{se:graphon}.
A graphon $Q$ is \emph{finitely forcible} if there are finitely many graphs $F_1,\ldots,F_\ell$ such that 
for every graphon $Q'$, if $p(F_i,Q)=p(F_i,Q')$ for all $i\in [\ell]$, then
$Q$ and $Q'$ are weakly isomorphic.

First, we need the following result which is Lemma~11 in~\cite{GlebovGrzesikKlimosovaKral15jctb} (except it is obtained by complementing all graphs and using our language of partite limits).

\begin{lemma}\label{lm:11} If $\OPT=\{\V x\}$ consists of a single element $\V x$ then there is $\ell_0$ (in fact, we can take $\ell_0=k$ where $k$ is as in the definition of $\lambda$) such that, for any $\V y\in\PC$ with $y_0=x_0$, if $p(\O{K_i}, \V x)=p(\O{K_i}, \V y)$ for every $2\leq i\leq \ell_0$ then $\V y=\V x$.
\end{lemma}

\bpf Our $\V x$ corresponds to a graphon $Q_{\V x}$. The fact that $\V x$ is the unique element of $\OPT$ is equivalent to saying that the equations $p(\O{P_3},Q)=0$ (the induced density of triples spanning exactly one edge) and $\lambda(Q)=\lambdamax$ force $Q$ to be $Q_{\V x}$ up to weak isomorphism in the space of graphons. In particular, $Q_{\V x}$ is finitely forcible. The constraint $p(\O{P_3},Q)=0$ forces $Q \in \C Q$ (that is, to be a complete partite graphon) and thus automatically forces $p(F,Q)=0$ for every graph $F$ which is not complete partite so we can ignore all such induced densities.

Thus the equation $\lambda(Q)=\lambdamax$ can be viewed as involving only induced densities of complete partite graphs on at most $k$ vertices. We claim that it can be equivalently rewritten as some polynomial in $x_0$ and induced densities of independent sets of size at most $k$. Then, supposing that the claim is true, if $Q_{\V y}\in\C Q$ has $y_0=x_0$ and the same induced densities of $\O{K_2},\ldots,\O{K_k}$ as $Q_{\V x}$, then $Q_{\V y}$ and $Q_{\V x}$ are weakly isomorphic and thus $\V y=\V x$.

It remains to prove the claim. For this, it suffices to prove that for any complete partite graph $F=K_{a_1,\ldots,a_\ell}$ with vertex set $[k]$ and with $\ell$ parts, for all $\V x \in \PC$, we have that $p(F,Q_{\V x})$ is some polynomial of $x_0$ and $p(\O{K_2},Q_{\V x}),\ldots,p(\O{K_k},Q_{\V x})$. 
The claim is clear for $\ell=1$ so assume $2 \leq \ell \leq k$.
Assume that $a_1,\ldots,a_\ell$ are in non-increasing order, and let $t \in [\ell]$ be the largest integer such that $a_1,\ldots,a_{t} \geq 2$.
Recall the analytic formula~(\ref{eq:lambdaan}) for $p(F,Q_{\V x})$.
We have
\begin{align}
\nonumber S_{a_1,\ldots,a_\ell}(\V x)
&= S_{a_1}(\V x) S_{a_2,\ldots,a_{\ell}}(\V x) - S_{a_1+a_2,a_3,\ldots,a_{\ell}}(\V x)-S_{a_2,a_1+a_3,\ldots,a_{\ell}}(\V x)-\ldots\\
\label{eq:Sind} &\hspace{1cm}\ldots - S_{a_2,\ldots,a_{\ell-1},a_1+a_{\ell}}(\V x)
\end{align}
and for every $a \geq 2$ we have $p(\overline{K_a},Q_{\V x})=S_a(\V x)$.
The claim now follows by induction on $\ell$. 
Indeed, every $S_{a_1,\ldots,a_{\ell-s}}(\V x)$ can be expressed as a polynomial of $S_a(\V x)$ for $2 \leq a \leq k$, by~(\ref{eq:Sind}) and induction, as required.
\epf

We need the following easy generalisation of Lemma~\ref{lm:11}.

\begin{lemma}\label{lm:11+} If $\OPT$ is finite then there is $\ell_0$ such that, for every $\V x\in\OPT $ and every $\V y\in\PC$ with $y_0=x_0$, if $\V x$ and $\V y$ have the same induced density of $\O{K_i}$ for every $1\leq i\leq \ell_0$ then $\V y=\V x$.
\end{lemma}

\bpf For every pair $\V z,\V z'\in\OPT$ there is some graph $F$ such that $p(F,\V z)\not=p(F,\V z')$. Indeed, since $\V z\not=\V z'$, their graphons $Q_{\V z}, Q_{\V z'}$ are not weakly isomorphic and thus have a different induced density of some graph $F$. Of course, this $F$ has to be complete partite (otherwise its induced density in both $\V z$ and $\V z'$ is zero). Let $F_1,\ldots,F_m$ be all such graphs $F$ where $m \leq \binom{|\OPT|}{2}$.
Let $\ell_0 := k+2\max_{i \in [m]} v(F_i)$.
Now let $\V x$ and $\V y$ be as in the lemma.

Consider the new optimisation problem where we maximise 
$$\lambda'(\V z):=\lambda(\V z)-\sum_{i=1}^m (p(F_i,\V z)-p(F_i,\V x))^2.
$$ Again, as in the proof of Lemma~\ref{lm:11}, $\lambda'$ can be written as a polynomial of $x_0$ and induced densities of anticliques on at most $\ell_0$ vertices. Also, clearly, $\V x$ is the unique element of $\OPT(\lambda')$. Apply Lemma~\ref{lm:11} to $\OPT(\lambda')=\{\V x\}$.\epf

\bpf[Proof of Lemma~\ref{lm:ff}] Let $\ell_0$ be as in Lemma~\ref{lm:11+}. It is enough to show that, for every $\V x\in\OPT$, there are at most $m:=\ell_0$ distinct non-zero values among $x_1,x_2,\ldots$ (then since $|\OPT|<\infty$ the lemma trivially follows).

Suppose on the contrary that $x_{i_1},\ldots,x_{i_{m+1}}$ are all positive and distinct for some $1 \leq i_1 < \ldots < i_{m+1}$.
Without loss of generality, assume that these are the smallest such indices we could have chosen.
Consider unknown variables $y_{i_1},\ldots,y_{i_{m+1}}$ and set $y_i:=x_i$ for every other $i\geq 1$. We get a contradiction to our choice of $\ell_0$ if we show that there is a choice of  $y_{i_1},\ldots,y_{i_{m+1}}>0$ such that
 \begin{equation}\label{eq:DMoment}
 \sum_{j=1}^{m+1} y_{i_j}^d=\sum_{j=1}^{m+1} x_{i_j}^d,\quad \mbox{for every $d=1,\ldots,m$},
 \end{equation}
 but the reordering $\V y'$ of $\V y$ (so that $y_1'\geq y_2'\geq \ldots$ and $y_0'=y_0$) is not equal to $\V x$. (Indeed, then $\V y'\in\PC$ by the case $d=1$ of~\eqref{eq:DMoment} and it satisfies $p(\O{K_d},\V y')=p(\O{K_d},\V x)$ for every $d=2,\ldots,\ell_0$ by the corresponding case of~\eqref{eq:DMoment}.)
 
Consider the map $g:\I R^m\times \I R\to \I R^m$ which sends $(z_1,\ldots,z_{m+1})$ to $(\sum_{j=1}^{m+1} z_j^d)_{d=1}^m$. The Jacobian of $g(\cdot,x_{i_{m+1}}):\I R^m\to\I R^m$, which sends
$\V z\in\I R^m$ to $g(\V z,x_{i_{m+1}})$, has non-zero determinant at $\V z_0:=(x_{i_1},\ldots,x_{i_m})$.
Indeed, the $(s,t)$-entry of the Jacobian at $(z_1,\ldots,z_m)$ is $s z_t^{s-1}$ so if we divide its $s$-th row by $s$ we obtain the Vandermonde matrix of $z_1,\ldots,z_m$, so its determinant is $m! \prod_{1\leq s<t\leq m} (z_s-z_t)$ which is non-zero at $\V z=\V z_0$.

Thus the Jacobian of $g(\cdot,x_{i_{m+1}})$ is invertible.
By the Implicit Function Theorem, for every choice of $y_{i_{m+1}}$ sufficiently close to $x_{i_{m+1}}$ 
there is a continuous choice of $(y_{i_1},\ldots,y_{i_m})$ close to $(x_{i_1},\ldots,x_{i_m})$ satisfying~\eqref{eq:DMoment}. 
Choose such a $y_{i_{m+1}}$ not equal to any $x_j$ and such that $y_{i_1},\ldots,y_{i_m}$ are all positive. Then the reordering $\V y'$ of the obtained sequence $\V y$ is not equal to $\V x$, giving the desired contradiction.
\epf

\section{The proof of Theorem~\ref{thm-main}}\label{sec-main}

In the first part of the proof, we find a suitable `hypothetical counterexample' $H$ on $h$ vertices (Claim~\ref{cl:findH}). This means that $H$ is very close to being optimal ($\lambda(H)$ is almost as large as $\lambda(h)$), but it is comparatively far from being complete partite (though it is important that $H$ is not \emph{too far} from being complete partite, and also that $H$ is very large).
Using~\Syone, given a candidate for $H$ which has too many imperfections, we can incrementally symmetrise it until this is no longer the case, and without decreasing $\lambda$.

In the second part of the proof (Claim~\ref{cl-vd}), we use the strictness of $\lambda$ to obtain a contradiction. We compare $H$ with the graph $H'$ obtained by removing all imperfections (roughly speaking $H'$ is the closest complete partite graph to $H$). The ratios of part sizes of $H'$ are necessarily close to some $\V x \in \OPT$.
The contradiction will come from the fact that $\lambda(H')-\lambda(H)$ is too large (which implies that $H$ is actually far from optimal).
We would like to argue that $\lambda(H)-\lambda(H')$ can be approximated looking at each \emph{wrong pair} $e \in W := E(H)\bigtriangleup E(H')$ separately and summing its contribution to the function.
This need not be true if $e$ is incident to many other wrong pairs, so instead we consider two families of wrong pairs: those incident to vertices in $B$, which are those with high degree in $W$, and the collection $E'$ of remaining wrong pairs.
The fact that each $e \in E'$ has a large contribution to $\lambda(H)-\lambda(H')$ will follow from~\Sone: namely that $\flip_{i_1i_2}\lambda(\V x)$ is large, where $i_1,i_2$ are the indices of the parts where $e$ lies.
The fact that the edges incident to each $v \in B$ have a large contribution to $\lambda(H)-\lambda(H')$ is slightly more involved.
For this we use~\Sytwo~to symmetrise the neighbourhood of $v$, and, depending on the attachment of $v$ in the resulting graph, the required conclusion will follow from~\Sone~(if it is `canonical') and~\Stwo~(otherwise).

The following lemma will be useful when comparing $\lambda$
evaluated on a complete partite graph with $\lambda$ evaluated on the same graph with a few imperfections.

\begin{lemma}\label{compare}
Let $c>0$ and let $\gamma : \mathcal{G}_k \to \mathbb{R}$ be fixed.
Let $H,H'$ be graphs on the same vertex set of size $h$, where $h$ is large and $H'\in \CP$ has $\PC$-structure $V_0,V_1,\ldots,V_m$.
Write $R := E(H)\bigtriangleup E(H')$ and given $x \in V(H')$, write $p(x)$ for the index of the part of $H'$ containing $x$. 
Define
$$
\xi_0 := k^2|R|c/h^2,\quad \xi_1 := 2\gammamax k^4|R|^2/h^4,\quad \xi_2 := 2\gammamax k^3|R|\Delta(R)/h^3.
$$
Then $\lambda(H')-\lambda(H)$ is
\begin{itemize}
\item[(i)] at least $\xi_0/2 -\xi_1-\xi_2$ if $\flip_{p(x)p(y)}\lambda(\V x_{H'}) \geq c$ for all $xy \in R$;
\item[(ii)] at least $\xi_0/2 -\xi_2$ if $\flip_{p(x)p(y)}\lambda(\V x_{H'}) \geq c$ for all $xy \in R$ and $R$ is a star;
\item[(iii)] at most $\xi_0 +\xi_1+\xi_2$ if $\flip_{p(x)p(y)}\lambda(\V x_{H'}) \leq c$ for all $xy \in R$.
\end{itemize}
\end{lemma}

\bpf
Write
$
S := \lambda(H')-\lambda(H) = \binom{h}{k}^{-1} \sum_{X \in \binom{V}{k}}\left(\gamma(H'[X])-\gamma(H[X])\right)
$
and
\begin{align*}
\nonumber S_1 &:= \binom{h}{k}^{-1}\sum_{xy \in R} \sum_{X \in \binom{V}{k}:\lbrace x,y\rbrace \subseteq X}(\gamma(H'[X])-\gamma((H' \oplus xy)[X]))\\
&= \frac{\binom{h-2}{k-2}}{\binom{h}{k}}\sum_{xy \in R}\flip_{p(x)p(y)}\lambda(H')
= \frac{\binom{k}{2}}{\binom{h}{2}}\sum_{xy \in R}\left(\flip_{p(x)p(y)}\lambda(\V x_{H'})+o(1)\right).
\end{align*}
Then
$
\binom{h}{k}|S-S_1| \leq \sum_{X \in I_1}2\gammamax$ where
$I_1 := \lbrace X \in \binom{V}{k}: |R \cap \binom{X}{2}| \geq 2\rbrace$.
The number of $X$ that contain two disjoint pairs from $R$ is
at most $|R|^2\cdot {h-4\choose k-4}$. The number of $X$ containing two adjacent pairs from $R$ is at most
$|R|\cdot \Delta(R)\cdot 
{h-3\choose k-3}$.
So
\begin{equation}\nonumber
\frac{|I_1|}{\binom{h}{k}} \leq \frac{|R|^2\binom{h-4}{k-4} + |R| \Delta(R)\binom{h-3}{k-3}}{\binom{h}{k}} \leq \frac{|R|^2k^4}{h^4} + \frac{|R| \Delta(R) k^3}{h^3}.
\end{equation}
All three parts follow immediately, noting for (ii)~that when $R$ is a star it has no disjoint pairs.
\epf

We now have all the tools in place to prove our main theorem.

\bpf[Proof of Theorem~\ref{thm-main}]
Let $\lambda$ be a symmetrisable graph parameter as in~(\ref{eq:gamma}). 
Note that $\lambda$ is not identically $0$ (otherwise $\OPT$ is infinite).
Lemma~\ref{lm:ff} implies that there exists $\beta>0$ such that $\OPT \subseteq \PC_\beta$.
So $|\supp(\V x)| \leq 1/\beta$ for all $\V x \in \OPT$.

Suppose that $\lambda$ is strict with parameter $c>0$. Without loss of generality we may assume that $c \ll \beta,1/\gammamax,1/k$.
We want to show that there exists a constant $C>0$ such that for every graph $G$ on at least $1/C$ vertices, there exists a complete partite graph $H$ on the same vertex set such that $\done(G,H)\leq C(\lambda(v(G))-\lambda(G))$.
Suppose that this is false.
That is, there exists a sequence of counterexamples $(G_n)_{n}$ with $v_n:=v(G_n)\to \infty$, such that
\begin{equation}\label{eq-GnCP}
1 \geq d_n:=\done(G_n,\CP_{v_n})>n(\lambda(v_n)-\lambda(G_n)),\quad\text{so}
\end{equation}
\begin{equation}\label{eq-Gnmax}
\lambda(v_n)\geq \lambda(G_n)> \lambda(v_n)-\frac{1}{n},
\end{equation}
and thus $ \lambda(G_n)-\lambda(v_n) \to 0$.

Using the graphs $G_n$, we now find a large graph $H$ which is almost optimal and has a small but \emph{comparatively large} number of imperfections.

\begin{claim}\label{cl:findH}
For all $\eps>0$, there exists $\eps'>0$ such that the following holds.
For all $N>0$, there exist $\V x \in \OPT$ and a graph $H$ on vertex set $[h]$ such that $h > N$, $\dedit(H,\V x) \leq 2\eps$ and $\lambda(H) \geq \lambda(h)-1/N$. Further,
$
\dedit(H,\CP_h) \geq \min\lbrace \eps', N(\lambda(h)-\lambda(H)) \rbrace
$.
\end{claim}

\bcpf
We consider two cases depending on whether $(d_n)_{n}$ contains a subsequence converging to $0$.
If it does not, then our counterexamples are eventually always far from being complete partite. In this case we perform an additional step of symmetrising each $G_n$ to obtain a graph which has a controlled number of imperfections; this number will be a small fraction of $v_n^2$.
In the other case, the counterexamples are becoming gradually more like complete partite graphs so the number of imperfections could be subquadratic (in $v_n$).

\medskip
\nib{Case 1}: $(d_n)_{n}$ does not contain a subsequence converging to $0$.

\medskip
\noindent
In this case, there exists $\xi>0$ such that $d_n \geq \xi$ for all sufficiently large $n$.
Since we are free to make $\eps$ and $\xi$ smaller we may assume without loss of generality that $\xi=\eps$.
Further, we may assume that $d_n \geq \eps$ for all $n \in \mathbb{N}$.

Let $V_n := V(G_n)$. Property \Syone~(applied with parameter $\eps$) implies that there exists $n_0=n_0(\eps)$ such that for each $n \geq n_0$, we can find a sequence $G_{n,0},G_{n,1},\ldots,G_{n,m(n)}$ of graphs on $V_n$ such that $G_{n,0} := G_n$; $G_n' := G_{n,m(n)}$ is complete partite; for all $i \in [m(n)]$, we have $\lambda(G_{n,i-1}) \leq \lambda(G_{n,i})$; and $\dedit(G_{n,i-1},G_{n,i}) \leq \done(G_{n,i-1},G_{n,i}) < \eps$.
By~\eqref{eq-Gnmax}, we have for all $0 \leq i \leq m(n)$ that $\lambda(v_n) \geq \lambda(G_{n,i}) \geq \lambda(G_n) > \lambda(v_n)-1/n$.

Let $\V y_n := \V x_{G_n'}$. By choosing a convergent subsequence since $(\PC,\dedit)$ is compact, we may assume that $\V y_n$ converges to some $\V y \in \PC$. But $\lambda(\V y_n)\rightarrow \lambdamax$, so $\V y \in \OPT$ by the continuity of $\lambda$.
By definition, $\dedit(G_{n,0},\OPT) \geq \dedit(G_{n,0},\CP_{v_n}) = d_n \geq \eps$ and
$\dedit(G_{n,m(n)},\OPT) \rightarrow \dedit(\V y,\OPT) = 0$.
Let $t$ be the largest element of $[m(n)]$ such that
$\dedit(G_{n,t},\OPT) \geq \eps$, and let $J_n := G_{n,t}$.
By increasing $n_0$, we can assume that $t<m(n)$.
Then
$$
\dedit(J_n,\OPT) \leq \dedit(G_{n,t},G_{n,t+1})+\dedit(G_{n,t+1},\OPT)  < 2\eps.
$$
That is, $\dedit(J_n,\OPT) \in [\eps,2\eps]$.
Let $\V x_n \in \OPT$ be such that $\dedit(J_n,\V x_n) = \dedit(J_n,\OPT)$.
 We claim that there exists $\eps'>0$ for which $p_n := \dedit(J_n,\CP_{v_n}) \geq \eps'$ for all sufficiently large $n$.

Indeed, if the claim is not true, then by passing to a subsequence we may assume that $p_n \rightarrow 0$.
For each $n$, pick a complete partite graph $P_n$ on $v_n$ vertices with $\dedit(J_{n},P_n) = p_n$.
Let $\V z_n := \V x_{P_n} \in \PC$ be the sequence that encodes the part ratios of $P_n$.
We can pass to a subsequence of $n$ such that $\V z_n$ converges to some $\V z \in \PC$;
then $\lambda(\V z) = \lim_{n\rightarrow \infty}\lambda(P_n) = \lambdamax$.
Thus $\V z \in \OPT$.
However, by Lemma~\ref{lm:edit},
\begin{align*}
\eps &\leq \dedit(J_{n},\OPT) \leq \dedit(J_{n},\V z) \leq \dedit(J_{n},P_n) + \dedit(\V z_n,\V z)\leq p_n + o(1) \rightarrow 0,
\end{align*}
a contradiction.

This $\eps'$ satisfies the lemma. Indeed, for any given $N>0$, choose $n>N$ sufficiently large so that $v_n > N$ and $\dedit(J_n,\V x_n) \in [\eps,2\eps]$ and $\dedit(J_n,\CP_{v_n}) \geq \eps'$.
Then we can set $\V x := \V x_n$ and $H := J_n$ and $h := v_n$, since $\lambda(J_n)\geq \lambda(v_n)-1/n\geq \lambda(h)-1/N$.
The claim is proved in this case.

\medskip
\nib{Case 2:}  $(d_n)_{n}$ contains a subsequence $(d_{n_i})_{i}$ such that $d_{n_i}\rightarrow 0$ as $i \rightarrow \infty$.

\medskip
\noindent
Assume without loss of generality that $(d_n)_{n} \rightarrow 0$. Therefore, there exists a sequence $(\V x_n)_{n}$ with $\V x_n \in \PC$ such that $\dedit(G_n,\V x_n)\rightarrow 0$. By choosing a convergent subsequence of $(\V x_n)_{n}$, we may assume that the sequence itself converges to  some $\V x \in \PC$. Then for sufficiently large $n$,
$\dedit(G_n,\V x)\leq \dedit(G_n,\V x_n)+\dedit(\V x_n,\V x) \rightarrow 0$.
Then the continuity  of $\lambda$ with respect to $\dedit$ and~\eqref{eq-Gnmax} imply that $\V x\in\OPT$.
We can choose $n$ sufficiently large so that, by~(\ref{eq-GnCP}), $H := G_n$ satisfies all the required properties in Claim~\ref{cl:findH} (where, for concreteness, we let $\eps':= 1$).
This completes the proof of the claim.
\ecpf

Choose an additional constant $0 < \eta \ll c$.
Obtain $\eps>0$ by applying Corollary~\ref{unicont} with $\eta^2,6\eps$ playing the roles of $\eps,\delta$ respectively.
We may assume that $\eps \ll \eta$.
Claim~\ref{cl:findH} furnishes us with an $\eps'>0$ which we may assume satisfies $\eps' \ll \eps$.
Now choose $N \in \mathbb{N}$ such that $1/N \ll \eps'$.
We have the following hierarchy:
\begin{equation}\label{hierarchy}
 0 <1/N \ll \eps' \ll \eps \ll \eta \ll  c \ll \beta,1/\gammamax,1/k.
\end{equation}
Apply Claim~\ref{cl:findH} to yield an $\V x \in \OPT$ and a graph $H$ on $h \geq N$ vertices.
Let us list some properties of $\V x$ (which will be all we need from now on):
\begin{itemize}
\item[\xsupp] $m := |\supp(\V x)| \leq 1/\beta$.
\item[\xbd] $\V x \in \PC_\beta$.
\item[\xone] $\flip_{i_1i_2}\lambda(\V x) \geq c$ for all $i_1,i_2 \in \supp^*(\V x)$.
\item[\xtwo] $\join_{b,\alpha}\lambda(\V x) \geq c((1-\alpha)x_0+ \min_{i \in \supp^*(\V x)}w_i)$, where
$$
\textstyle w_i = \mathds{1}_{i > 0}b_ix_i + \sum_{j\in\supp^*(\V x)\setminus \lbrace 0,i\rbrace}(1-b_j)x_j,
$$ 
for all $b:\mathbb{N}\to\{0,1\}$ and $\alpha \in [0,1]$.
\item[\xcts] Whenever $\V y \in \PC$ satisfies $\dedit(\V x, \V y) \leq 6\eps$ and $\supp^*(\V y)=\supp^*(\V x) =: S$, we have that $|f(\V x)-f(\V y)| \leq \eta^2$ for all choices of $i_1,i_2 \in S$, $b:\mathbb{N}\rightarrow\lbrace 0,1\rbrace$, $\alpha\in[0,1]$ and $f \in \lbrace \lambda,\flip_{i_1i_2}\lambda,\lambda(\cdot,(b,\alpha)),\join_{b,\alpha}\lambda\rbrace$.
\item[\xH] $\dedit(\V x, H) \leq 2\eps$.
\end{itemize}
Properties~\xsupp~and~\xbd~follow immediately from $|\OPT|<\infty$ and Lemma~\ref{lm:ff}. Properties~\xone~and~\xtwo~follow since $\lambda$ is strict with parameter $c$.
Property~\xcts~follows from our choice of $\eps$ and the fact from Corollary~\ref{unicont} that any $f$ in this family of functions is uniformly equicontinuous. Property~\xH~is a direct consequence of Claim~\ref{cl:findH}.

Let $\mathcal{H}$ be the family of $h$-vertex graphs with $\PC$-structure $(V_i: i \in \supp^*(\V x))$, that is, $V_0$ (if it exists) is a clique, $V_i$ is a non-empty independent set for all $i \in [m] = \supp^*(\V x)\setminus \{0\}$ and $(V_i,V_j)$ is complete for every distinct $i,j \in \supp^*(\V x)$.

Among all graphs in $\mathcal{H}$, let $H'$ be one whose edit distance $\dedit$ to $H$ is minimised, with $\PC$-structure $(V_i: i \in \supp^*(\V x))$ as above, where $V(H')=[h]=\bigcup\{V_i : i \in \supp^*(\V x)\}$. Define $W:=([h],E(H)\bigtriangleup E(H'))$, and call the edges of $W$ \emph{wrong}. By the definition of $H'$,~\xH~and Lemma~\ref{lm:edit}, we have that $\dedit(H,H') \leq \dedit(H,\V x)+O(1/h)\leq 3\eps$. Consequently $e(W)=\Done(H,H')\leq 3h^2/2 \cdot \dedit(H,H') \leq 5\eps h^2$. 
Let $\V v$ be the vector of part ratios in $H'$,~i.e.\ $\V v := \V x_{H'}$.

Then
\begin{equation}\label{eq-vx}
	\dedit(\V v,\V x)=\dedit(H',\V x)\le\dedit(H,H')+\dedit(H,\V x)\stackrel{\xH}{\le} 5\eps.
\end{equation}
Note that, by~\xbd,~this implies $v_i = |V_i|/h \geq c/2$ for all $i \in [m]$.
Call a vertex $x$ \emph{bad} if it is incident to at least $\eta h$ wrong pairs, i.e.\ $d_{W}(x)\ge\eta h$.
Let $B$ consist of all bad vertices and $B^c=[h]\setminus B$ of all \emph{good} (i.e.\ not bad) vertices.
Let also $E' := E(W[B^c])$ and $e' := |E'|$.
By definition of $B$ and that $\eps\ll \eta$, we have
\begin{equation}\label{eq-B}
e' \leq e(W) \leq 5\eps h^2\quad\text{and}\quad	|B|\le\frac{2e(W)}{\eta h}\leq \frac{10\eps}{\eta}\cdot h\le\sqrt{\eps} h.
\end{equation}
For a vertex $v$ of $H'$ let $H'\oplus v$ denote the graph obtained from $H'$ by removing every edge containing $v$ and then for all $y \in [h]\setminus \lbrace v \rbrace$ adding the edge $vy$ if and only if $y \in N_H(v)$. The heart of the proof is the following claim.

\begin{claim}\label{cl-vd}
   The following statements hold.
\begin{itemize}
		\item[(i)] $\dedit(H',H)\leq 2\left(\frac{|B|}{h}+\frac{e'}{h^2}\right)$.
\item[(ii)] For every $v\in B$, $\lambda(H')-\lambda(H'\oplus v)\geq \frac{kc\eta^{3/2}}{3h}$.

		\item[(iii)] $\lambda(H')-\lambda(H)\geq \eta^2\left(\frac{|B|}{h}+\frac{e'}{h^2}\right)$.
	\end{itemize}
\end{claim}
We first see how this claim completes the proof of Theorem~\ref{thm-main}.
We have by Claim~\ref{cl:findH} that
\begin{eqnarray}
\nonumber	\frac{1}{N} &\geq& \lambda(h)-\lambda(H)\ge\lambda(H')-\lambda(H)\stackrel{(iii)}{\ge} \eta^2\left(\frac{|B|}{h}+\frac{e'}{h^2}\right)\stackrel{(i)}{\ge}\frac{\eta^2}{2}\dedit(H',H) \geq \frac{\eta^2}{2}\dedit(H,\CP_h)\\
\label{eq-contra} &\geq& \frac{\eta^2}{2}\min\lbrace \eps', N(\lambda(h)-\lambda(H)) \rbrace.
\end{eqnarray}
If $\eps' \leq N(\lambda(h)-\lambda(H))$, then considering the first and last terms of~(\ref{eq-contra}) gives $1/N \geq \eta^2\eps'/2$, a contradiction to our choice of $N$ (i.e.\ (\ref{hierarchy})).
If instead $\eps' > N(\lambda(h)-\lambda(H))$, then considering the second and last terms of~(\ref{eq-contra}) gives $1 \geq \eta^2N/2$, also a contradiction to our choice of $N$.
Thus Theorem~\ref{thm-main} holds given Claim~\ref{cl-vd}.

\bpf[Proof of Claim~\ref{cl-vd}.]
For~(i), we see that
\begin{eqnarray*}
	\dedit(H',H)\leq \done(H',H) \leq \frac{\sum_{v\in B}d_W(v)+e'}{h^2/2}\leq \frac{2|B|}{h}+\frac{2e'}{h^2}.
\end{eqnarray*}

For (ii), fix an arbitrary $v\in B$, and let $p(v) \in \supp^*(\V x)$ be such that $v \in V_{p(v)}$. Let $\C H$ consist of all graphs $G$ on $[h] = V(H')$ with $G-v=H'-v$ and for each $i\in [m]$, either $N_G(v)\supseteq V_i\setminus\{v\}$ or $N_G(v)\cap (V_i\setminus\{v\})=\varnothing$ (and with arbitrary attachment to $V_0$). 
That is, either $v$ is adjacent to every vertex or no vertices in each part $V_1,\ldots,V_m$.
For brevity, let $H'_v:=H'\oplus v$.

Apply \Sytwo~with parameter $\eps$ to $H'_v$ at $v$ to obtain a sequence of graphs $H'_v=: H_0,H_1,\ldots,H_r\in\C H$, such that $H_i-v = H'_v-v$ for all $i \in [r]$; $\lambda(H_{i-1}) \leq \lambda(H_i)$; and $\Done^v(H_{i-1},H_i) \leq \eps(h-1)$ for all $i \in [r]$,
where here for any two graphs $J,J'$ which differ only at a vertex $v$, we define $\Done^v(J,J')$ to be the minimum number of edits of pairs containing $v$ to make $J$ equal to $J'$.

By the definition of $\mathcal{H}$, there are $b : [m]\rightarrow \lbrace 0,1\rbrace$ and $0 \leq \alpha \leq 1$ such that $b(i)=1$ if $N_{H_r}(v) \supseteq V_i\setminus\lbrace v \rbrace$ and $b(i)=0$ if $N_{H_r}(v) \cap (V_i\setminus \lbrace v \rbrace)=\varnothing$; and $d_{H_r}(v,V_0) = \lfloor\alpha|V_0|\rfloor$ (if $V_0=\varnothing$ we let $\alpha := 1$).
We consider two cases depending on $(b,\alpha)$: in Case~1 the attachment of $v$ in $H_r$ is very different to any vertex in $H'-v$, and in Case~2 it is similar.

\medskip
\nib{Case 1:} At least one of the following holds: (a) $x_0>0$ and $\alpha < 1-\eta/2$; (b) $|b^{-1}(0)|\geq 2$; (c) $x_0=0$ and $b^{-1}(0) = \varnothing$.

\medskip
\noindent
We will first show that
\begin{equation}\label{nablabalpha}
\frac{h}{k}\left(\lambda(H')-\lambda(H_r)\right) \geq \join_{b,\alpha}\lambda(\V x)-3\eta^2.
\end{equation}
For this, let $\V y$ be the vector of part ratios of $H'_r-v=H'-v$, i.e.\ $y_i = |V_i|/(h-1)$ if $i \in \lbrace 0,\ldots,m\rbrace\setminus\lbrace p(v)\rbrace$ and $y_{p(v)} = (|V_{p(v)}|-1)/(h-1)$. 
Then $H_r = (H'-v) +_{b,\alpha} v = G_{h-1,\V y} +_{b,\alpha} v$ and so
\begin{align*}
\frac{h}{k}\left(\lambda(H') - \lambda(H_r)\right) &= \frac{h}{k}\cdot\binom{h}{k}^{-1}\sum_{X \in \binom{V}{k} : X \ni v}\left(\gamma(H'[X]) - \gamma((G_{h-1,\V y}+_{b,\alpha}v)[X])\right)\\
&= \frac{h}{k}\cdot\frac{\binom{h-1}{k-1}}{\binom{h}{k}}\left( \lambda(H',v) - \lambda(G_{h-1,\V y}+_{b,\alpha}v,v)\right)\\
&= \lambda(\V v,(e_{p(v)},1)) - \lambda(\V y,(b,\alpha)) + O(1/h).
\end{align*}
Now $\dedit(\V x, \V y),\dedit(\V x, \V v) \leq 5\eps$. So we have
\begin{align*}
\lambda(\V v,(e_{p(v)},1)) - \lambda(\V y,(b,\alpha)) &\stackrel{\xcts}{\geq} \lambda(\V x,(e_{p(v)},1)) - \lambda(\V x,(b,\alpha)) - 2\eta^2 = \join_{b,\alpha}\lambda(\V x)- 2\eta^2,
\end{align*}
where the final equality follows from Proposition~\ref{lagmult}(i).
This proves~(\ref{nablabalpha}).
Now,
$$
\lambda(H')-\lambda(H'\oplus v) \geq \lambda(H') - \lambda(H_r) \stackrel{(\ref{nablabalpha})}{\geq} \frac{k}{h}\left(\join_{b,\alpha}\lambda(\V x)-3\eta^2\right)
$$
so to complete the proof of the claim, it suffices to show that $\join_{b,\alpha}\lambda(\V x) \geq c\eta^{3/2}/2$.

We will use the lower bound on $\join_{b,\alpha}\lambda(\V x)$ given by~\xtwo, and that $x_i \geq \beta > c$ for all $i \in \supp^*(\V x)$ from~\xbd.
Suppose first that~(a) holds.
Since each term in the expression for $w_i$ is non-negative, we have for all $i \in \supp^*(\V x)=\lbrace 0,\ldots,m\rbrace$ that
$\join_{b,\alpha}\lambda(\V x) \geq c(1-\alpha)x_0 \geq c\beta\eta/2> c\eta^{3/2}/2$, as required.
Suppose secondly that~(b) holds.
Then for every $i \in [m]$, either $b_i=1$ or $b_j=0$ for some $j \in [m]\setminus\lbrace i \rbrace$. So $w_i \geq \min_{j \in [m]}x_j \geq \beta$ for all $0 \leq i \leq m$, as required.
Finally, if~(c) holds, $x_0=0$, $b=(1,1,\ldots)$,
$\supp^*(\V x) = [m]$ and $w_i = x_i \geq \beta$ for all $i \in [m]$, as required.

\medskip
\nib{Case 2:} Either (a) $x_0=0$ and $|b^{-1}(0)|=1$, or (b) $x_0>0$, $\alpha > 1-\eta/2$, and $|b^{-1}(0)|\leq 1$.

\medskip
\noindent
Notice that Cases~1 and~2 are the only possible outcomes (recalling that if $x_0=0$, then $\alpha=1$).
For all $0 \leq i \leq r$, let
$$
d_i := \min_{j \in \supp^*(\V x)}\Done^v(H_i, (H'-v) +_{e_j,1} v),
$$
i.e.\ the smallest number of edits at $v$ needed to move $v$ into some part in $H_i$.
Now, $d_0 = d_W(v) \geq \eta h$.
On the other hand, $d_r$ is comparatively small: if~(a) holds, then $d_r=0$, and if~(b) holds, then $d_r \leq \eta h/2$.
So we can choose the largest integer $0 \leq t < r$ such that $d_t \geq \eta h$, and let $H^* := H_t$ and $d^* := d_t$.
Let $j^* \in \supp^*(\V x)$ be such that $d^* = \Done^v(H_t, \widetilde{H})$ where $\widetilde{H} := (H'-v) +_{e_{j^*},1} v$.
So
$
\lambda(H^*) \geq \lambda(H' \oplus v)
$
and additionally
$$
\eta h \leq d^* \leq \Done^v(H_t, H_{t+1}) + d_{t+1} \leq \eps(h-1)+\eta h \leq 2\eta h.
$$
So one must make between $\eta h$ and $2\eta h$ edits to $H^*$ at $v$ to move $v$ to the $j^*$-th part, and the (complete partite) graph obtained in this way is $\widetilde{H}$.
Since $\dedit(\widetilde{H},H')\cdot h^2/2 \leq \Done(\widetilde{H},H') \leq h$, we have by~(\ref{eq-vx}) that $\dedit(\V x_{\widetilde{H}},\V x) \leq \dedit(\widetilde{H},H') + \dedit(H',\V x) \leq 5\eps + 2/h \leq 6\eps$.
Now, \xcts~and \xone~imply that
for each of the $d^*$ vertices $u$ for which $uv$ was flipped, we have
$\flip_{p(u)p(v)}\lambda(\V x_{\widetilde{H}}) \geq \flip_{p(u)p(v)}\lambda(\V x)-\eta^2 \geq c/2$.
Lemma~\ref{compare}(ii) implies that
$
\lambda(\widetilde{H})-\lambda(H^*) \geq k^2d^*c/2h^2 - 2\gammamax k^3 (d^*)^2/h^3
$.
So
$$
\lambda(H')-\lambda(H' \oplus v)  \geq \lambda(\widetilde{H})- \lambda(H^*) + O(1/h) \geq k^2\eta^{3/2} c/3h,
$$
as required for (ii).

\medskip
For~(iii), our task is to obtain a suitable lower bound on
$T:=\lambda(H')-\lambda(H)$.
Notice that the only $k$-sets $X$ contributing to $T$ are those containing some $e\in W$.
Let
	\begin{eqnarray*}
		T_0&:=\sum_{v\in B}(\lambda(H')-\lambda(H'\oplus v))\quad\text{and}\quad T' := \lambda(H')-\lambda(H'\bigtriangleup E').
\end{eqnarray*}
In a similar fashion to part~(ii), we will first give lower bounds for $T_0,T'$ respectively, and then show that $T$ is well-approximated by $T_0+T'$.
First
consider $T_0$. By Claim~\ref{cl-vd}(ii), we have
$T_0\geq |B|k\eta^{3/2} c/(3h)$.
Now consider $T'$.
Again $\flip_{p(x)p(y)}\lambda(\V x_{H'}) \geq c/2$ for all $xy \in E'$, so
Lemma~\ref{compare} and~(\ref{eq-B}) imply that
\begin{align*}
T' &\geq k^2e'/h^2 \cdot (c/2 - 50\gammamax k^2 \eps^2 - 2\gammamax k\eta) \geq \frac{k^2ce'}{4h^2}.
\end{align*}
For the final step, note that
$\binom{h}{k}|T - T_0 - T'|  \leq \sum_{X \in I_0}2\gammamax$,
where
$$
 I_0 = \left\lbrace X \in \binom{V}{k} : |X \cap B| \geq 2 \text{ or } |X \cap B|,e(W[X\setminus B])\geq 1\right\rbrace.
$$
But
\begin{equation}\label{I0}
\frac{|I_0|}{\binom{h}{k}} \leq \frac{|B|^2\binom{h-2}{k-2} + |B| e'\binom{h-3}{k-3}}{\binom{h}{k}} \stackrel{(\ref{eq-B})}{\leq} \frac{|B|}{h}\left(2k^2\sqrt{\eps} + 5k^3\eps\right) \leq \eps^{1/3}T_0.
\end{equation}
Thus
$$
T \geq T_0 + T' - \frac{2\gammamax}{\binom{h}{k}}|I_0| \stackrel{(\ref{I0})}{\geq} T_0/2 +T' \geq \frac{|B|k\eta^{3/2} c}{6h} + \frac{k^2ce'}{4h^2} \geq \eta^2\left(\frac{|B|}{h}+\frac{e'}{h^2}\right),
$$
as desired.
This completes the proof of Claim~\ref{cl-vd}.
\ecpf

Thus we complete the proof of Theorem~\ref{thm-main}.
\epf


\section{Applications to inducibility}\label{sec-app}

First we prove Lemma~\ref{lem-symm} which is essentially Theorem~1 in~\cite{SchelpThomason98}.

\bpf[Proof of Lemma~\ref{lem-symm}] In fact, we can require that $|E(G_{i-1})\bigtriangleup E(G_i)|$ is at most $n-1$ (resp.\ at most $1$) in \Syone~(resp. \Sytwo) for every graph $G$ of every order $n\geq k$.

Let us show~\Syone. Initially, let $H:=G$ and let $\C V=\{V_1,\ldots,V_n\}$ be the partition of $V(H)$ into singletons. At each stage, every part of $\C V$ will consist of twin vertices, i.e.\ vertices with identical neighbourhoods (in particular, every part is an independent set). We will modify the current graph $H$ and the current partition $\C V=\{V_1,\ldots,V_s\}$ so that at each step $\lambda$ does not decrease while the affected edges are incident to a single vertex.

If for each $1\leq i<j\leq s$, $H[V_i,V_j]$ is complete bipartite, then $H$ is a complete partite graph so we stop. Otherwise, pick $i<j$, $x\in V_i$ and $y\in V_j$ such that $xy\not\in E(H)$. Let $X=N_H(x)$ and $Y=N_H(y)$. Fix a complete partite graph $F$. Note that every $A \subseteq V$ with $H[A] \cong F$ is one of the following four kinds: (1) $x \in A$, $y \notin A$; (2) $x \notin A$, $y \in A$; (3) $x \in A$, $y \in A$; and (4) $x \notin A$, $y \notin A$. Given $H-x-y$, we can thus write 
\begin{equation}\nonumber
p(F,H)=f_F(X)+f_F(Y)+g_F(X\cap Y, V\setminus(X\cup Y))+C_F
\end{equation}
for some constant $C_F>0$ and functions $f_F$ and $g_F$. Here $f_F(X)$ (resp. $f_F(Y)$) counts the number of copies of $F$ of type (1) (resp. type (2)) as this depends only on $X$ (resp. $Y$).
For disjoint $U,W$, we define $g_F(U,W)$ to be the number of copies of any graph $J$ with $V(J) \subseteq U \cup W$ such that by adding two new vertices $z,z'$ to $J$ and adding edges $\{uz,uz': u \in U\}$ to $J$ we obtain a copy of $F$.
Observe that if $\{x,y\}\cup V(J)$ induces a copy of $F$ in $H$ as above, then $x$ and $y$ are in the same partite set, $U\cap V(J) \subseteq X \cap Y$ and $W \cap V(J)\cap (X \cup Y)=\varnothing$.
Thus $g_F(X\cap Y, V\setminus(X\cup Y))$ counts type (3) copies. The type (4) count is a constant depending only on $H-x-y$.
Then, letting $f=\sum_{F}c_F\cdot f_F$ and defining $g,C$ similarly, we have
\begin{equation}\label{eq-4kinds}
\lambda(H)=f(X)+f(Y)+g(X\cap Y, V\setminus(X\cup Y))+C.
\end{equation}
Notice that $g(\cdot,\cdot)$ is non-decreasing in both arguments, that is,
\begin{equation}\label{eq-mono}
g(U,W)\leq g(U',W'), \quad\forall~U\subseteq U', W\subseteq W'.
\end{equation}
Indeed, if $F$ is a clique, then no copy of $F$ contains both $x$ and $y$, and $c_F\geq 0$ otherwise. 

Suppose that $f(X)\geq f(Y)$, let
$H_{xy}$ be the graph obtained from $H$ by making $y$ a clone of~$x$. Let $H'=H_{xy}$ and let $\C V'$ be obtained from $\C V$ by moving $y$ to the part containing $x$. It satisfies all the claimed properties as 
\begin{eqnarray*}
\lambda(H_{xy})&=&2f(X)+g(X,V\setminus X)+C\\
&\overset{(\ref{eq-4kinds}),(\ref{eq-mono})}{\ge}&
f(X)+f(Y)+g(X\cap Y, V\setminus(X\cup Y))+C=\lambda(H).
\end{eqnarray*}

Finally, it remains to argue that one can avoid infinite cycles. The rule for breaking ties $f(X)=f(Y)$ with e.g.~$|V_i|\geq |V_j|$ is to take $H'=H_{xy}$. This strictly increases $\sum_{V\in\C V} |V|^2\in [n,n^3]$ so that are at most $n^3$ steps where $\lambda$ stays constant. (In fact, one can bound the total number of steps by $1+2+\ldots+n-1={n\choose 2}$: if there are currently $i\geq 2$ groups and we merge one group entirely into another, then we can do this by moving at most $n-i+1$ vertices.)

Let us show \Sytwo. Given $G$ and $z$ as in the property, we have a partition consisting of all partite sets in $G-z$ and $z$ will always stay a single part. Given any partite set $V_i$ of $G-z$, we can partition vertices $V_i=V_i'\cup V_i''$ depending on their adjacency to $z$, say $V_i'\subseteq N(z)$. 
Start with this initial partition into parts $V_i'$ and $V_i''$.
Fix arbitrary non-adjacent vertices $x\in V_i', y\in V_i''$, note that~\eqref{eq-4kinds} and~\eqref{eq-mono} still hold. 
If $f(X) > f(Y)$, take $G'=G_{xy}$. If $f(X) < f(Y)$, take $G'=G_{yx}$. The rule for breaking ties is again to clone the vertex from the larger part: if $f(X)=f(Y)$ and, say, $|V_i'| \geq |V_i''|$, take $G'=G_{xy}$.
Otherwise, take $G'=G_{yx}$. Note that $G'$ differs from $G$ only in one pair.
As before, $\lambda$ has not decreased.
Then redefine $V_i', V_i''$ and repeat the process.
The final graph has $N(z) \subseteq V_i$ or $N(z) \cap V_i=\varnothing$.
(Note that each tie $f(X)=f(Y)$ strictly increases $|V_i'|^2+|V_i''|^2$ so as before there are at most $n^3$ steps where $\lambda$ stays constant, so there are no infinite cycles.)
Repeating this for all $i$, we make at most $n$ steps in total, and the resulting graph is as desired.
\epf

\subsection{Proofs of Theorems~\ref{thm-Kst}--\ref{thm-K311}}

Since by Lemma~\ref{lem-symm}, $p(F,\cdot)$ is symmetrisable whenever $F$ is complete partite,
to prove Theorems~\ref{thm-Kst}--\ref{thm-K311} it suffices to determine $\OPT$ (if it is not already known),
and then check that $p(F,\cdot)$ is strict.
The result then follows from Theorem~\ref{thm-main}.

In all cases, $\OPT$ consists of a single point, and checking strictness is generally straightforward (it is slightly more involved for $F=K_{1,t}$).
However, determining $\OPT$ where it is not already known, for $F=K_{2,1,1,1}$ and $F=K_{3,1,1}$, is challenging
and we are required to solve a polynomial optimisation problem.
We use computer-assisted semidefinite programming to solve the last problem.

\bpf[Proof of Theorem~\ref{thm-Kst}]
Assume $s \leq t$.
First we collect  some facts about the function $f_{s,t}$ defined on $[0,1]$ given by $f_{s,t}(\alpha)=\alpha^s(1-\alpha)^t+\alpha^t(1-\alpha)^s$, recalling that for $s<t$, $f_{s,t}(\alpha)= \binom{t+s}{t}^{-1}p(K_{s,t},(\alpha,1-\alpha,0,\ldots))$ for $\alpha \in [\frac{1}{2},1]$, and $f_{s,s}$ can similarly be expressed with a factor of $\frac{1}{2}$ on the right-hand side:
\begin{itemize}
\item[(i)] If $s \geq \binom{t-s}{2}$ then the unique maximum of $f_{s,t}$ in $[\frac{1}{2},1]$ is $\frac{1}{2}$.
\item[(ii)] If $s< \binom{t-s}{2}$, then $f_{s,t}'$ has a single root in $(\frac{1}{2},1)$, which corresponds to a maximum, has $\frac{1}{2}$ a root corresponding to a minimum, and has no other roots in $[\frac12,1)$.
\item[(iii)] If $\alpha\in [\frac12,1]$ maximises $f_{1,t}$, then $1-\alpha>\frac{1}{t+1}$.
\item[(iv)] $\max_{\alpha \in [0,1]}(t+1)\alpha^t(1-\alpha) = \left(\frac{t}{t+1}\right)^t$, attained uniquely at $\frac{t}{t+1}$.
\end{itemize}
Note that $f_{s,t}$ is symmetric about $\alpha=\frac{1}{2}$. For~(i) and~(ii),
we just follow the proof of~\cite[Theorem~3]{BrownSidorenko94}. We have
$$
f_{s,t}'(\alpha) = \alpha^{t}(1-\alpha)^{s-1}h\left(\frac{1-\alpha}{\alpha}\right)\quad\text{where }
h(x)= sx^{t-s+1}-tx^{t-s}+tx-s.
$$
Assume first that $s \geq \binom{t-s}{2}$.
Setting $x=1+\eps$ for $\eps \geq 0$, one can show that $h(x)>0$, so $f_{s,t}$ is non-decreasing in $[0,\frac{1}{2}]$, and thus the unique maximum of $f_{s,t}$ in $[0,1]$ is at $\frac{1}{2}$, as required for (i).
(In the calculation in~\cite[Theorem~3]{BrownSidorenko94}, it is shown that $h(x) \geq 0$, but there is equality in the first inequality only if $t-s=1$, but in this case the final inequality is strict.) Note that~\cite{BrownSidorenko94} uses $(t,s+t)$ and $(a,b)$ instead of our $(s,t)$.
Assume secondly that $s<\binom{t-s}{2}$.
Following the remarks after~\cite[Theorem~5]{BrownSidorenko94}, it suffices to show that $h$ has a single root in $(0,1)$. This is a consequence of $h(0)<0$, $h(1)=0$, $h'(1)<0$ and $h''(x)<0$ for all $x \in (0,1)$, as required for~(ii).

Next we show that~(iii) holds. If $t=2,3$ (i.e.\ $1 \geq \binom{t-1}{2}$), then (i) implies that $1-\alpha=\frac{1}{2}$, as required.
If $t \geq 4$, then by~(ii), $f_{1,t}'$ has a unique root (i.e.\ $1-\alpha$) in $(0,\frac{1}{2})$, corresponding to a maximum and $\frac{1}{2}$ is a root corresponding to a minimum.
Thus $f_{1,t}'(x)>0$ for $x \in (0,1-\alpha)$ and $f_{1,t}'(x)<0$ for $x \in (1-\alpha,\frac{1}{2})$.
One can check that $f'(\frac{1}{t+1})>0$, which gives $1-\alpha > \frac{1}{t+1}$.
This proves~(iii). Property~(iv) can be easily checked via differentiation:
indeed, $\frac{d}{d\alpha}\alpha^t(1-\alpha) = \alpha^{t-1}(t-(t+1)\alpha)$ is strictly positive for $0<\alpha < \frac{t}{t+1}$, equals $0$ at $\alpha=\frac{t}{t+1}$, and is strictly negative for $\alpha>\frac{t}{t+1}$.

Now we show that $\OPT = \{(\alpha,1-\alpha,0,\ldots)\}$ with $\alpha \geq \frac{1}{2}$
(where $f_{s,t}(\alpha)$ is uniquely maximised).
This was essentially proved by Brown and Sidorenko~\cite{BrownSidorenko94}. They do not prove the uniqueness of the optimal element but this can be extracted from their proof, so we only give a sketch of how to do this here.

First we claim that if $G$ is a complete multipartite graph on $n$ vertices 
whose two largest parts of sizes $n_r,n_{r-1}$ satisfy $n_r,n_{r-1}=\Omega(n)$ and $n-n_r-n_{r-1}=\Omega(n)$, then by merging parts, we increase the number of induced copies of $K_{s,t}$ by $\Omega(n^{s+t})$. 
Indeed, to see the claim, fix $\eps>0$ and suppose $G=K_{n_1,n_2,\ldots,n_r}$ with $r\geq 3$ and $n_1\leq n_2\le\ldots\leq n_r$ with $\sum_{i\in[r]}n_i=n$, $n_{r-1} \geq \eps n$, $\sum_{i\in[r-2]}n_i \geq \eps n$, and $G'=K_{n_1+n_2,n_3,\ldots,n_r}$.
It is shown in~\cite[Proposition~2]{BrownSidorenko94} that merging the two smallest parts in any complete multipartite graph with at least three parts does not decrease the number of induced copies of $K_{s,t}$.
Thus in $G$ we can successively merge two smallest parts until we obtain a graph $G''$ with exactly three parts, of sizes $m_1 \leq m_2 \leq m_3$ with $m_1 \geq \eps n$.
Now merge the parts of size $m_1$ and $m_2$ to obtain a complete bipartite graph $G'''$.
Then
\begin{eqnarray*}
&&I(K_{s,t},G''')-I(K_{s,t},G)\\
&\geq& I(K_{s,t},G''')-I(K_{s,t},G'')\\
&=&{m_1+m_2\choose s}{m_3\choose t}+{m_1+m_2\choose t}{m_3\choose s}\\
	&&\vspace{2cm}-{m_1\choose s}{m_3\choose t}-{m_1\choose t}{m_3\choose s}-{m_2\choose s}{m_3\choose t}\\
	&&\vspace{2cm}-{m_2\choose t}{m_3\choose s}-{m_1\choose s}{m_2\choose t}-{m_1\choose t}{m_2\choose s}\\
	&\ge&{m_1+m_2\choose s}{m_2\choose t}+{m_1+m_2\choose t}{m_2\choose s}-{m_1\choose s}{m_2\choose t}-{m_1\choose t}{m_2\choose s}\\
	&&\vspace{2cm}-{m_2\choose s}{m_2\choose t}-{m_2\choose t}{m_2\choose s}-{m_1\choose s}{m_2\choose t}-{m_1\choose t}{m_2\choose s}\\
	&=&\frac{1}{t!s!}\left( m_2^t((m_1+m_2)^s-m_1^s-m_2^s-m_1^s)+m_2^s((m_1+m_2)^t-m_1^t-m_2^t-m_1^t)\right)+O(n^{s+t-1}).
\end{eqnarray*}
To prove the claim, it suffices to prove that this is at least $O(\eps)n^{s+t}$ for all $(s,t)$.
Neglecting the $O(n^{s+t-1})$ error terms, the quantity in the last line is at least
$$
\begin{cases}
			\frac{1}{t!s!}\left(m_1^s(t-2)m_1^t + m_1^t(s-2)m_1^s\right) \geq \frac{1}{t!s!}m_1^{s+t} & \text{if }s+t \geq 5\\
            \frac{1}{3!}\left(2m_1m_2^3+3m_1^2m_2^2-m_1^3\right) \geq \frac{2}{3}m_1^4 & \text{if }(s,t)=(1,3)\\
	\frac{1}{2!2!}2m_2^2(2m_1m_2-m_1^2) \geq \frac{1}{2}m_1^4 & \text{if }(s,t)=(2,2)\\
	\frac{1}{2}(m_2^2m_1-m_1^2m_2)  \geq (m_2-m_1) m_1^2 & \text{if }(s,t)=(1,2).
		 \end{cases}
$$
Since $m_1 \geq \eps n$, this proves the claim unless $(s,t)=(1,2)$ and $m_2-m_1 < \eps n$.
In this case, we have $I(K_{1,2},G''')-I(K_{1,2},G'') = m_1^2m_3-m_1^3+O(\eps)n^3=\mu^2(1-3\mu)n^3+O(\eps)n^3$ where $\mu := m_1/n$. We are done if $\mu < \frac{1}{3}-\eps$. If not, $G''$ has three parts of size $\frac{1}{3}\pm \eps$ which is far from optimal, by comparing to the complete balanced bipartite graph.
This completes the proof of the claim.

Thus, if $s \geq 2$, then every $\V x \in \OPT$ has exactly two non-zero entries which sum to $1$, as required.
We want to show this also holds for $s=1$.
For this, we only need to show that there is no $0<x \leq 1$ for which $\V x = (x,0,\ldots)$ is optimal. Indeed, $\lambda(\V x)=p(K_{1,t},\V x)=(t+1)x^t(1-x) \leq \left(\frac{t}{t+1}\right)^t$ by~(iv).
But
$
\lambda(K_{1,t},(\textstyle{\frac{t}{t+1}},\textstyle{\frac{1}{t+1}},0,\ldots))=(t+1)f_{1,t}(\frac{t}{t+1}) = 2\left(\frac{t}{t+1}\right)^t
$,
so $\V x \not\in \OPT$.

We have shown that every element of $\OPT$ is of the form $(\alpha,1-\alpha,0,\ldots)$ for some $\alpha \in [\frac{1}{2},1]$.
By~(i) and~(ii), $f_{s,t}$ has a unique maximum in $[\frac{1}{2},1]$.
Thus $\OPT$ contains a unique element $\V x := (\alpha,\beta,0,\ldots)$ (where from now on we write $\beta := 1-\alpha$).
It remains to show that there is $c=c(\lambda)>0$ such that~\Sone~and~\Stwo~hold, where $\lambda(\cdot):= p(K_{s,t},\cdot)$.

Let $G:=G_{n,\V x}$.
First we check~\Sone.
A non-edge between two partite sets is not contained in any induced copy of $K_{s,t}$, nor is an edge within a partite set.
So $\flip_{xy}\lambda(G)$ is the number of copies of $F$ in $G$ containing the pair $xy$ divided by $\binom{n-2}{s+t-2}$.
Rather roughly, this is always at least $\beta^{s+t-2}+o(1)$ as $n\to\infty$.

Now we check~\Stwo.~We have $x_0 = 0$ and $\supp^*(\V x) = \lbrace 1,2\rbrace$. Since $x_0=0$, given any $(b(1),b(2)) =: (b_1,b_2)\in\lbrace 0,1\rbrace^2$, we are required to show that
$$
\join_{b,1}\lambda(\V x) =   \lambda(\V x)- \lim_{n\rightarrow\infty}\lambda(G_{n,\V x}+_{b} u,u) \geq c\min\lbrace  b_1\alpha + (1-b_2)\beta, b_2\beta+(1-b_1)\alpha\rbrace.
$$
Recall that as usual the right-hand side equals $0$ if $(b_1,b_2)\in\lbrace (0,1),(1,0)\rbrace$, so the inequality is trivially true.
If $(b_1,b_2)=(0,0)$, then $u$ lies in no copies of $K_{s,t}$ in $G_{n,\V x}+_{b,1} u$;
similarly if $(b_1,b_2)=(1,1)$ and $s \geq 2$.
So we may assume that $(b_1,b_2,s)=(1,1,1)$, and we need to show $\join_{b,1}\lambda(\V x) \geq c\beta$ (recalling $\alpha \geq \beta$).
We have
\begin{align*}
\lambda(\V x) &= (t+1)(\alpha\beta^t + \alpha^t\beta), \quad \lambda(\V x, (b,1)) = \alpha^t+\beta^t.
\end{align*}
Recall that $\frac{\partial \lambda(\V x)}{\partial x_i}=(t+1)\lambda(\V x)$ for $i=1,2$ by Proposition~\ref{lagmult}(ii). We have
\begin{align*}
\frac{\partial \lambda(\V x)}{\partial x_1} &= (t+1)(\beta^t + t\alpha^{t-1}\beta) \quad \text{and}\quad \frac{\partial \lambda(\V x)}{\partial x_2} = (t+1)(t\alpha\beta^{t-1} + \alpha^t), \text{ and so}
\end{align*}
\begin{align*}
2\join_{b,1}\lambda(\V x) &= 2\lambda(\V x) - 2\lambda(\V x,(b,1)) = 
\frac{1}{t+1}\left(\frac{\partial \lambda(\V x)}{\partial x_1}+\frac{\partial \lambda(\V x)}{\partial x_2}\right)-2(\alpha^t+\beta^t)\\
&=\alpha^{t-1}(t\beta-\alpha)+\beta^{t-1}(t\alpha-\beta)
=\alpha^{t-1}((t+1)\beta-1) + \beta^{t-1}((t+1)\alpha-1).
\end{align*}
It suffices to show that $\beta > \frac{1}{t+1}$,
since then writing $\eps=\beta-\frac{1}{t+1}$, we have
$\join_{b,1}\lambda(\V x) \geq (t+1)^{1-t}\eps$.
This follows from~(iii),
completing the proof. 
\epf

\bpf[Proof of Theorem~\ref{thm-Krt}]
First we show that for $F:=K_r(t)$ with $t>1+\log r$, we have $\OPT = \{\V x\}$ where $\V x = (\frac{1}{r},\ldots,\frac{1}{r},0,\ldots)$ and $x_0=0$.
This essentially follows from~\cite{BrownSidorenko94} where it is proved that $\V x$ lies in $\OPT$ (but without proving uniqueness), and~\cite[Theorem~13]{BollobasEgawaHarrisJin95} where it is proved that the Tur\'an graph with $r$ parts is the unique extremal graph; but as in the proof of Theorem~\ref{thm-Kst} we again need to make some modifications.
Write $\lambda(\cdot):=p(K_r(t),\cdot)$, and observe that $\lambda(\V y)=(tr)!/(r!(t!)^r) \cdot S_{t,\ldots,t}(\V y)$ where $t$ is repeated $r$ times.
The method of Lagrange multipliers,~\cite[Proposition~7]{BrownSidorenko94} shows that every $\V y$ which maximises $\lambda$ has exactly $r$ non-zero entries (which sum to $1$, since $t>1$).
Thus it suffices to show that $S^t_r(x_1,\ldots,x_r) := x_1^t\ldots x_r^t$ over all $x_1,\ldots,x_r > 0$ with $x_1+\ldots +x_r= 1$ is uniquely maximised by $(\frac{1}{r},\ldots,\frac{1}{r})$.
This is easy to see for $r=2$. Suppose it is not true for some $r \geq  3$, and so $x_1 \neq x_2$, say.
Then $S^t_r(x_1,\ldots,x_r)>S^t_r(\frac{x_1+x_2}{2},\frac{x_1+x_2}{2},x_3,\ldots,x_r)$, a contradiction.

We have proved that $\OPT=\{ \V x\}$.
It remains to show that there is $c=c(\lambda)>0$ such that~\Sone~and~\Stwo~hold.

Note that \Sone~is immediate as a non-edge between two partite sets is not contained in any induced copy of $F$, and an edge within a partite set is not contained in any induced copy of~$F$.
As in the proof of Theorem~\ref{thm-Kst}, this means that every $\flip_{ij}\lambda(\V x)=\lambda(\V x)$, which is always at least $\left(\frac{1}{r}\right)^{tr-2}+o(1)$.
Similarly \Stwo~is immediate since any vertex in $G_{n,\V x}$ without neighbours in at least two parts does not have a $K_{r-1}$ in its neighbourhood so does not lie in any copies of $K_r(t)$,
and any dominating vertex clearly lies in no copies of $K_r(t)$.
So again $\join_{b,1}\lambda(\V x) = \lambda(\V x)$
whenever $b \neq e_i$ for any $i \in [r]$.
\epf

\bpf[Proof of Theorem~\ref{thm-K2111}] 
Let us show that $\lambdamax$ is $\lambda_0:=\frac{525}{1024}$ and the vector $\V a = (\frac{1}{8},\ldots,\frac{1}{8},0,\ldots)$, which is the limit of $K^8_{n/8,\ldots,n/8}$, is the unique maximiser. (Here, $K^\ell_{n_1,\ldots,n_\ell}$ is the complete $\ell$-partite graph with parts of size $n_1,\ldots,n_\ell$.)

Let $\V x\in\OPT$ be arbitrary. At some places it will be convenient to use the language of finite graphs. So let $n$ be large and let $G=G_{n,\V x}$ be a realisation of $\V x$ with $\PC$-structure $V_0,\ldots,V_m$.

Let us show that there are $\ell \in \mathbb{N}$, $p:=(1-x_0)/\ell$ and 
a sequence $\V x = \V x_0,\ldots, \V x_t = (p,\ldots,p,0,\ldots) \in \OPT$ where for all $j \in [t]$, $x_{j,0}=x_0$, the entries of $\V x_j$ are obtained by replacing some non-zero $x_{j-1,i_1},x_{j-1,i_2}$ in the entries of $\V x_{j-1}$ with $z_j := x_{j-1,i_1}+x_{j-1,i_2}$, and any $\V y_j$ obtained by replacing $x_{j-1,i_1},x_{j-1,i_2}$ in the entries of $\V x_{j-1}$ with non-negative reals that sum to $z_j$ is also in $\OPT$.



Indeed, suppose there are non-zero $x_i \neq x_j$, and let $i+j$ be minimal with this property. Each copy of $F=K_{2,1,1,1}$ intersects each $V_i\cup V_j$ in at most 3 vertices. Thus, if we fix the rest of $\V x$, fix $s=x_i+x_j$ and
vary $a=x_i/s$ between $0$ and $1$, then
the number of copies of $F$ is given by a polynomial $p(a)$ of degree at most $2$ which is symmetric around $\frac{1}{2}$: $p(a)=p(1-a)$. (Note that the number of copies of $F$ having $2+1$ vertices in $V_i\cup V_j$ is a constant times $a^2(1-a)+(1-a)^2a=a-a^2$, which has no $a^3$ term, i.e.\ is also a quadratic polynomial.)
If $p$ is not constant, then by symmetry, it follows that $p'(\frac{1}{2})=0$ and since $p'$ is a linear function of $a$, this is the only root. Thus $p$ is maximised at $0$, $1$ or $\frac{1}{2}$ and we can strictly increase $p$, a contradiction.
Thus $p$ is constant, and any $\V z$ obtained from $\V x$ by replacing $x_i$ and $x_j$ by one or two new entries whose sum of sizes is $x_i+x_j$ is in $\OPT$ (corresponding to taking any value $0 \leq a \leq 1$).
We let $i_1=i$ and $i_2=j$, giving $\V x_1$. Then $\V x_1$ and any $\V y_1$ as described lie in $\OPT$.
If we cannot take $t=1$, then $\V x_1$ has unequal non-zero entries and we can repeat the above.
It remains to check that this process terminates.
If not, since we can always merge the largest part with the next (non-equal non-zero) largest part, for all $\eps>0$ there is some $m = m(\eps)>0$ such that $x_{m,1}>1-x_0-\eps$ (recalling $x_{m,0}=x_0$).
Then $p(F,\V x_m) = 5!(1-x_0)^2/2 \cdot x_0^3/6 + O(\eps)$ which is maximised when $x_0=\frac{3}{5}$, with value $\frac{216}{625}< \frac{525}{1024}$, a contradiction.

Let $\V y = \V x_t$, so $\V y$ has $\ell$ equal non-clique parts each of ratio $p=(1-y)/\ell$ for some $y\in [0,1]$. Thus $p(F,\V y)$ is equal to $h_\ell(y)$, where
 \begin{eqnarray*}
 h_\ell(y)&:=&5!\, \ell\ \frac{p^2}2\left( \frac{(1-p)^3}{3!}-(\ell-1)\frac{p^2}2(1-2p)-(\ell-1)\frac{p^3}{3!}\right).
 \end{eqnarray*}
 Indeed, we first choose one part $V_i$ where two non-adjacent vertices go ($\ell {pn\choose 2}$ choices). Then the other 3 vertices of $F$ have to go outside of $V_i$ (${(1-p)n\choose 3}$ choices) except we have to rule out two (exclusive) cases: exactly two of them are in some $V_j$ ($(\ell-1){pn\choose 2}(n-2pn)$ choices) and all three of them are in some $V_j$ ($(\ell-1){pn\choose 3}$ choices). 
We have $h_\ell'(y)=\frac{10}{\ell^4}(y-1)q(y)$ where
$$
q(y)=-30+49\ell-21\ell^2+2 \ell^3+90 y-123 \ell y+33 \ell^2 y-90 y^2+99 \ell y^2-12 \ell^2 y^2+30 y^3-25 \ell y^3.
$$
 
We claim that for each $\ell\geq 8$, the function $h_\ell$ is strictly monotone decreasing (that is, the optimal $y$ is 0 meaning that the clique part is empty). 
So it suffices to show that $q(y)>0$ for $\ell \geq 8$ and $y \in [0,1]$.
We have that $q$ is positive at its endpoints: $q(0)=42+97(\ell-8)+27(\ell-8)^2+2(\ell-8)^3$ and $q(1)=2\ell^3$.
So if $q(y)<0$ for some $y \in [0,1]$, then $q'$ has a root in $[0,1]$.
However, the quadratic polynomial $q'$ has a negative coefficient at $y^2$ and is positive at endpoints: $q'(0)=1218+405(\ell-8)+33(\ell-8)^2$ and $q'(1)=9\ell^2$,
so there is no such root.
(These symbolic calculations can be found in \texttt{2111.nb} in the ancillary folder.)

We claim that $k \in \mathbb{ R}[\ell]$ given by $k(\ell)=h_\ell(0)$ is decreasing for $\ell \geq 8$.
That is, out of all $\V y$ with at least $8$ non-zero entries which are all equal, the unique extremal $\V y$ is $\V a$. Indeed, $k'(\ell)= -10j(\ell)/\ell^5$ where $j(\ell) = (\ell-9)^3 + 15 (\ell-9)^2 + 60 (\ell-9) + 30$ so $k'(\ell)$ is decreasing for all $\ell \geq 9$, and also $k(9)=\frac{1120}{2187} < \frac{525}{1024} = k(8)$.
(See \texttt{2111.nb}.)

Let us show that none of $\ell\in [7]$ is optimal. Fix such an $\ell$. Direct calculations show that $h_\ell(0) < \lambda_0$ (while $y=1$ gives $K_n$ which has zero density of $F$). So it remains to investigate critical points, that is, $y \in (0,1)$ such that $h_\ell$ has derivative zero at $y$. Thus $q(y)=0$.

Introduce a new variable $z$ and define $p_1(y):=h_\ell'(y)$ and $p_2(y,z):=z-h_\ell(y)$. Thus if $y$ is a critical point with $\lambdamax=h_\ell(y)$ and we define $z:=h_\ell(y)$ then $(y,z)$ belongs to the variety $V=V(I)\subseteq \I R^2$ defined by the ideal $I:=\langle p_1,p_2\rangle$ generated by the polynomials $p_1,p_2$. By applying Buchberger's algorithm to $I$ (where we eliminate the variable $y$) we see that $J$, the intersection of $I$ with $\I R[z]$ (the set of polynomials that depend on $z$ only) is generated by one polynomial $q_\ell$, explicitly computed in \texttt{2111.nb} for every $\ell \in [7]$. We actually need only a part of the above claim, namely that there are polynomials $f_1,f_2\in \I R[y,z]$ such that we have a polynomial identity $q_\ell(z)=f_1(y,z)p_1(y)+f_2(y,z)p_2(y,z)$, that is, all terms on the right-hand side depending on $y$ cancel each other.


We have $q_1(z)=z(625z-216)$ which has roots at $z=0,\frac{216}{625}$, and $h_1(z_0)<\lambdamax$ for both roots $z_0$.
For each $2\leq l\leq 7$, the polynomial $q_\ell$ on inspection has the following properties: we have $q_\ell(z)=zr_\ell(z)$ where $r_\ell$ has degree at most $3$, the coefficient of the leading term of $r_\ell$ is positive, and furthermore, we have $r_\ell(0)\geq 0$, $r_\ell(\lambda_0)<0$ and $r_\ell(1)<0$. This implies that $r_\ell$ has no roots in $[\lambda_0,1]$ and hence $q_\ell$ has no roots in $(0,1]$. That is, it is impossible to have $\lambdamax>\lambda_0$ (because, as a graph density, $\lambdamax$ is at most 1). Thus $\lambdamax=\lambda_0$ and none of the polynomials can achieve $\lambdamax$ except when $\ell=8$ (with $\V y = (\frac{1}{8},\ldots,\frac{1}{8},0,\ldots)$ being the unique maximiser among $\V y \in \PC$ with $y_i=y_j$ for all $i,j \in [8]$).

If $\V y = \V x_t \neq \V x_0$, then $\V x_{t-1}$ exists, and it is of the form $\V x_{t-1} := (\frac{1}{8},\ldots,\frac{1}{8},a,\frac{1}{8}-a,0,\ldots)$ for some $\frac{1}{16} \leq a < \frac{1}{8}$, where $\frac{1}{8}$ is repeated $7$ times, and moreover the element of $\PC$ obtained by setting $a=\frac{1}{16}$, say, lies in $\OPT$.
A routine calculation shows this to be a contradiction (see \texttt{2111.nb}).

 Thus $\OPT=\{\V a\}$, where $\V a := (\frac{1}{8},\ldots,\frac{1}{8},0,\ldots)$ with $a_0=0$.

Finally, it remains to check strictness.
Let us check~\Sone.
First let $x,y$ be in different parts $V_i,V_j$ of $G=G_{n,\V a}$. Write $p=\frac{1}{8}$ and assume $8 | n$. 
Consider copies of $K_{2,1,1,1}$ that contain both $x$ and $y$, with $A$ denoting the $2$-element part.
Then in 
$G$, the edge $xy$ can be such that $x$ lies in $A$ ($pn \cdot\binom{6}{2}\cdot (pn)^2$) choices), $y$ lies in $A$ (the same), or neither $x$ nor $y$ lie in $A$ ($6\cdot\binom{pn}{2}\cdot 5pn$ choices).
In $G \oplus xy$, $xy$ is a non-edge and so we can only have $\{x,y\}$ playing the role of $A$ in $F$, so the number of such copies of $F$ is $\binom{6}{3}\left(pn\right)^3$. Thus for distinct $i,j \in [8]$ we have
$$
\flip_{ij}\lambda(\V a) = 3!p^3\left(2\binom{6}{2} + 15 - \binom{6}{3}\right) = \frac{150}{512}.
$$
Now let $x,y$ be in the same part $V_i$ of $G$. Then in $G$, the non-edge $xy$ lies in $\binom{7}{3}(pn)^3$ copies of~$F$.
In $G \oplus xy$, the edge $xy$ lies in $7\binom{pn}{2} \cdot 6pn$ copies of $F$. So
$$
\flip_{ii}\lambda(\V a) = 3!p^3\left(\binom{7}{3} - \frac{7\cdot 6}{2}\right) = \frac{84}{512},
$$
as required.

For~\Stwo, let $b: [8] \to \{0,1\}$ be such that $|\supp(b)|=k$. Then
$$
\lambda(\V a, (b,1)) =\lim_{n\to\infty}\binom{n}{4}^{-1}\left(\frac{n}{8}\right)^4 \left( (8-k)\binom{k}{3} + \binom{k}{2} (k-2)\cdot\frac{1}{2}\right) = \frac{4!}{8^4}\cdot\binom{k}{3}\left(\frac{19}{2}-k\right).
$$
Indeed, counting induced copies of $F$ in $G_{n,\V a} +_{b,1} u$ containing $u$: if $u$ plays the role of a vertex in $A$, then  we choose the other vertex from this set among any of the $8-k$ parts not adjacent to $u$, and then choose three distinct parts of the $k$ adjacent to $u$ to contain the other vertices.
If $u$ plays the role of a singleton, we choose two among $k$ parts for the other two singletons, and another for $A$ (dividing by two for both orders).
Routine calculations show that this is uniquely maximised (with value $\lambda_0$) when $k=7$, as required.
\epf

\bpf[Proof of Theorem~\ref{thm-K311}]
We will show that $\OPT=\{(\frac{3}{5},0,\ldots)\}$ and $\lambdamax = \frac{216}{625}$.
Let $G$ be a complete partite graph on $n$ vertices which maximises the number of $F:=K_{3,1,1}$.
Comparing $G$ to $(\frac{3}{5},0,0,\ldots)$, we have $P(F,G) \geq \frac{216}{625}\binom{n}{5} + O(n^4)$.
Suppose that $Y,Z$ are the two largest parts of $G$, with $|Y|=yn$, $|Z|=zn$ and $y \geq z$.
Let $S:=V(G)\setminus(Y \cup Z)$.

First, let us derive a contradiction from assuming that $z \geq \frac{2}{5}$.
Let $s:=1-y-z$, so $s \leq \frac{1}{5}$.
The number of copies of $F$ with at least three vertices in $S$ is at most
$$
n^5\left(\frac{s^5}{5!} + \frac{s^4(1-s)}{4!}+\frac{s^3}{3!}\cdot yz\right) \leq
n^5\left(\frac{s^5}{5!} + \frac{s^4}{4!}+\frac{s^3}{4!}\right) \leq n^5\frac{151}{600\cdot 625}.
$$
The number of copies of $F$ with exactly two vertices in $S$ is at most
\begin{equation}\label{eq:2S}
\binom{sn}{2}\left(\binom{yn}{3}+\binom{n-sn-yn}{3}\right) = \frac{n^5s^2}{12}(y^3+(1-s-y)^3) + O(n^4).
\end{equation}
We have $y \leq 1-s-\frac{2}{5}$ (since $z \geq \frac{2}{5}$) and $y>1-s-y$ (since $y \geq z$). For fixed $s$, the expression $y^3+(1-s-y)^3$ is maximised when $y$ is as large as possible.
Indeed, $y^3+(1-s-y)^3$ for $y \in \mathbb{R}$ is a quadratic polynomial whose coefficient at $y^2$ is positive and whose minimum is at $\frac{1-s}{2}$, and we have $y>\frac{1-s}{2}$.
So the expression in~(\ref{eq:2S}) is at most
$$
r(s)n^5 + O(n^4)\quad\text{where}\quad r(s) = \frac{s^2}{12}\left(\left(1-s-\frac{2}{5}\right)^3+\left(\frac{2}{5}\right)^3\right).
$$
We claim that $r'$ has no roots in $(0,\frac{1}{5}]$, which implies that $r(s)$ attains its maximum at $s=\frac{1}{5}$, of value $\frac{160}{600 \cdot 625}$.
Indeed, $r'(s)=\frac{s}{300}t(s)$ where $t(s)=-125s^3+180s^2-81s+14$.
Furthermore, $t(1)<0<t(\frac{4}{5})$ so $t$ has at least one root in $[\frac{4}{5},1]$. If the claim does not hold, then $t$ has three real roots, which are interlaced by the roots of the quadratic $t'(s)=-3(5s-3)(25s-9)$. The smallest root of $t'$ is $\frac{9}{25}>\frac{1}{5}$, and the coefficient of $s^3$ in $t$ is negative, so $t$ has a root in $(0,\frac{1}{5}]$ only if $t(\frac{9}{25})<0$, a contradiction.

Every other copy of $F$ has exactly $4$ vertices in $Y \cup Z$. So, writing $q := \frac{1-s}{y}$, their number is
$$
\binom{yn}{3}zn \cdot sn + \binom{zn}{3}\cdot yn \cdot sn = 
\frac{n^5}{6}(q^3(1-q)+(1-q)^3q)s(1-s)^4 + O(n^4)
$$
which, for $s \in [\frac{4}{5},1]$, is maximised when $(s,q)=(\frac{4}{5},\frac{1}{2})$, with value $\frac{640}{600}\cdot \frac{n^5}{625}+O(n^4)$.
So when $z \geq \frac{2}{5}$, we have $p(F,G) \leq \frac{5!}{600\cdot 625}(151+160+640) < \lambda_0$, and
we obtain the desired contradiction.

Assume from now on that $z < 2/5$.
Fix $v \in Z$. Let $p(F,G,v)$ be the number of copies of $F$ containing $v$. Then
\begin{align*}
P(F,G,v) &\leq p(v) := \binom{zn-1}{2}\left(\binom{(1-z)n}{2}-\binom{yn}{2}\right)+\binom{yn}{3}(1-y-z)n\\
&\hspace{2cm}+\frac{1}{3}\sum_{w \in S}\binom{n-1-d(w)}{2}(d(w)-zn).
\end{align*}
We would like a good upper bound for the last term.
Since $Z$ is the second largest part, we have that $(1-z)n \leq d(w) \leq n$ for all $w \in S$.
Now $f(x)=\frac{1}{2}(1-x)^2(x-z)$ is maximised when $x=x_0 := \frac{1}{3}(1+2z)$ and is decreasing on the interval $[x_0,1]$.
Since $z \leq \frac{2}{5}$, we have $x_0 \leq 1-z$, so $f$ defined in the range $[1-z,1]$ is maximised at $x=1-z$.
So the last term divided by $n^4$ is at most
$$
\frac{1}{3}\sum_{w \in S}f(d(w))n^{-1} + O(1/n) = \frac{1}{3}(1-y-z) f(1-z) + O(1/n).
$$
Define
\begin{align*}
h(y,z) &:=12\left( \frac{z^2}{4}((1-z)^2-y^2) + \frac{y^3}{6}(1-y-z) + \frac{1}{3}(1-y-z)f(1-z)-\frac{9}{625}\right)\\
&= 2y^3-2 y^4-2 y^3z+5 z^2-2 y z^2-3 y^2z^2-12 z^3+4 y z^3+7 z^4-\frac{108}{625}.
\end{align*}
By the above, $h(y,z) \geq 12(p(v)n^{-4} + O(1/n) - \frac{9}{625}) \geq O(1/n)$, that is, $h(y,z) \geq 0$ for all $0 \leq z \leq y$ with $z+y \leq 1$ and $z \leq \frac{2}{5}$. 
Let

$$
R := \{ (y,z) \in [0,1]^2: y \geq z, y+z \leq 1\}.
$$

\begin{claim}\label{cl:h}
For every $(y,z) \in R$ with $h(y,z) \geq 0$, we have that $y \geq \frac{3}{5}$.
\end{claim}

Suppose that the claim holds. 
Since $G$ is optimal, Proposition~\ref{lagmult} implies that $v$ has optimal attachment in $G$; that is, $P(F,G,v) = \binom{n-1}{4}\lambda(G,v) = \binom{n-1}{4}\lambdamax + O(n^3) \geq \frac{9}{625}n^4 + O(n^3)$.
Thus $h(y,z) \geq 0$ for the $y,z$ corresponding to $Y,Z$, since, as we have shown, $z \leq \frac{2}{5}$.
So $y \geq \frac{3}{5}$.
Consider the graph $H$ obtained by replacing $Z$ by a clique. 
Then we lose every copy of $F$ containing the $3$-independent set in $Z$ (and lose no other copies)
while we gain copies of $F$ with the $3$-independent set in $Y$ and the two other vertices in $Z$. So
\begin{align*}
\frac{P(F,G)-P(F,H)}{120n^5} &\leq\frac{z^3}{3!}\frac{(1-z)^2}{2} - \frac{(3/5)^3}{3!}\frac{z^2}{2}+o(1)
=\frac{z^2}{2}\left(\frac{z(1-z)^2}{6}-\frac{3^3}{6\cdot 5^3}\right)+o(1)\\
&\leq \frac{z^2}{2}\left(\frac{4}{6\cdot 3^3}-\frac{3^3}{6\cdot 5^3}\right)+o(1) \leq \frac{-229}{40500}z^2+o(1).
\end{align*}
This is a contradiction to the optimality of $G$ if $z=\Omega(1)$.
Thus $z=o(1)$ and, up to $o(n^2)$ edits, $G$ consists of an independent set of size $yn$ and $(1-y)n$ universal vertices.
So $p(F,G)=120(\frac{y^3}{3!}\frac{(1-y)^2}{2})+o(1)$. Ignoring the error term, this is uniquely maximised when $y=\frac{3}{5}$, with value $\frac{216}{625}$.
Then $\OPT=\{\V a\}$, where $\V a = (\frac{3}{5},0,\ldots)$.

So, in order to determine $\OPT$, it remains to prove Claim~\ref{cl:h}.

\bpf[Proof of Claim~\ref{cl:h}.]
First we consider $(y,z)$ on the boundary of $R$. If $z=0$ then $h(y,0)=\frac{-108}{625}+1250y^3(1-y)$
which is uniquely maximised when $y=\frac{3}{5}$.
If $y=z$ then $h(y,y)=y^2(2y-1)(2y-5)-\frac{108}{625}$ which is negative for $y \in [0,1]$.

Now we consider $(y,z)$ in the interior of $R$.
Let $(y_0,z_0)$ in the interior of $R$ be such that $h(y_0,z_0) \geq 0$ and $y_0$ is minimal with this property
(such a $y_0$ exists by compactness of $R$ and continuity of $h$).
Since $(y_0,z_0)$ is in the interior of $R$, we have $h(y_0,z_0)=0$ and $\frac{\partial h}{\partial z}(y_0,z_0)=0$
(otherwise we can find $z' \approx z_0$ with $h(y_0,z')>h(y_0,z_0)=0$ and by the continuity of $h$, $y' < y_0$ and $h(y',z') \geq 0$, contradicting the minimality of $y_0$).
Applying Buchberger's algorithm to eliminate $z$, we obtain a degree-$12$ polynomial $q$
such that $y$ satisfies $h(y,z)=0=\frac{\partial h}{\partial z}(y,z)$ only if $q(y)=0$ (see \texttt{311.nb}):
\begin{align*}
q(y) &:= -2500858044
+14506020000y
-18911610000y^2
-85830803750y^3
+545884288750y^4\\
&-1430659375000y^5
+4001212109375y^6
-12503827343750y^7
+30477566015625y^8\\
&-54597656250000y^9
+64171142578125y^{10}
-42002929687500y^{11}
+12102539062500y^{12}.
\end{align*}
Let $\alpha := \frac{272}{1000}$ and $R' := \{(y,z) \in (0,1]^2: z \leq y \leq \alpha\}$.
We claim that $p(y) := q(y+\alpha)$ is a positive polynomial.
Then $q(y) > 0$ for all $y \in R\setminus R'$, and hence $(y_0,z_0) \in R'$. 
For this, it suffices to show that there are polynomials $r_1(y), r_2(y)$ with non-negative coefficients 
satisfying $p(y)r_1(y) = r_2(y)$.
Once one fixes the degree $d$ of $r_1$, this amounts to solving a linear program, where $a_k$ is the $k$-th coefficient of $p$ and $b_k$ is the $k$-th (unknown) coefficient of $r_1$:
\begin{equation*}
\begin{array}{ll@{}ll}
\text{minimise}  & \displaystyle\sum\limits_{0 \leq k \leq d} b_{k} &\\
\text{subject to}& \displaystyle\sum_{\substack{j+k=i: \\0 \leq j \leq 12;\\ 0 \leq k \leq d}}a_jb_k > 0,   & \ \ \ i=0,1 ,\ldots, d+12,\\
                 &                                                b_{k} > 0, &\ \ \ k=0,1,\ldots, d.
\end{array}
\end{equation*}
In fact we only need a feasible solution, not an optimal one, so the objective function can be anything.
For degrees $d=1,2,\ldots$ we attempted this (using python) until we obtained a numerical solution for $d=16$. The following degree-16 polynomial was obtained by multiplying this solution by a fairly large power of 10 and rounding.
\begin{align*}
r_1(y) &= 405631585336 x^{16}+291048000156 x^{15}+172228102580 x^{14}+76577243592 x^{13}\\
&+32501733953 x^{12}+13576227809 x^{11}+5344727909 x^{10}+1954537506 x^9+737097269 x^8\\
&+264696828 x^7+90984085 x^6+30184081 x^5+10472958 x^4+3090485 x^3+1000538 x^2\\
&+206609 x+108298.
\end{align*}
Clearly its coefficients are positive and one can check (see \texttt{311.nb}) that the degree-28 polynomial $p(y)r_1(y)$ also has positive coefficients, as required.

Suppose we can find non-negative polynomials $s_0,\ldots,s_3$ in $y,z$ and positive $t \in \mathbb{Q}$
such that
$$
-h(y,z)-t-zs_1 - (y-z)s_2 - (\alpha-y)s_3=s_0,
$$
where a polynomial $p \in \mathbb{R}[y,z]$ is non-negative if $p(y,z) \geq 0$ whenever $y,z \geq 0$.
Then $-h(y,z) > 0$ on $R'$.
This will complete the proof of the claim.
Let $\underline{x}:=(1,y,z,y^2,yz,z^2)^\intercal$. To ensure that the $s_i$ are non-negative, it suffices to find positive semidefinite $6 \times 6$ matrices $Q_i$ such that $s_i(y,z)=\underline{x}^\intercal Q_i \underline{x}$.
For this, a sum-of-squares solver (we used the YALMIP Matlab toolbox~\cite{Lofberg04,Lofberg09} with SeDuMi~\cite{Sturm98}) numerically maximises $t$ such that the above equality holds; that is, we obtain $t' \approx 0.02$ and real matrices $Q_0',\ldots,Q_3'$ such that $-h(y,z)-t'-zs_1' - (y-z)s_2' - (\alpha-y)s_3' \approx s_0'$, where $s_i' = \underline{x}^\intercal Q_i' \underline{x}$.
Now let $Q_i$ be a (symmetric) rational approximation to $Q_i'$ for $i \in [3]$ and let $R_0$ be a rational approximation to $Q_0'$.
We obtain

\begin{align*}
R_0 &= \left(
\begin{array}{cccccc}
 \frac{47560627}{605583685} & -\frac{27288737}{128683162} & -\frac{5823553}{403766228} & -\frac{22660833}{166625377} & \frac{64761638}{445638833} & -\frac{10092851}{42370543} \\
 -\frac{27288737}{128683162} & \frac{412450960}{208083677} & -\frac{154126052}{222170865} & -\frac{123333398}{74059181} & -\frac{45208772}{76054353} & \frac{29997552}{77062243} \\
 -\frac{5823553}{403766228} & -\frac{154126052}{222170865} & \frac{56961038}{76246587} & \frac{75134651}{68479911} & -\frac{114623437}{74768701} & \frac{68436686}{157424595} \\
 -\frac{22660833}{166625377} & -\frac{123333398}{74059181} & \frac{75134651}{68479911} & \frac{231222579}{42911653} & -\frac{33046138}{90840815} & -\frac{27557233}{25108228} \\
 \frac{64761638}{445638833} & -\frac{45208772}{76054353} & -\frac{114623437}{74768701} & -\frac{33046138}{90840815} & \frac{142375474}{17195129} & -\frac{204334483}{99244906} \\
 -\frac{10092851}{42370543} & \frac{29997552}{77062243} & \frac{68436686}{157424595} & -\frac{27557233}{25108228} & -\frac{204334483}{99244906} & \frac{152251273}{45491357} \\
\end{array}
\right) \succeq 0\\
Q_1 &= \left(
\begin{array}{cccccc}
 \frac{113823133}{103564772} & -\frac{153720698}{116964597} & -\frac{514694857}{175951034} & -\frac{26958123}{134065612} & -\frac{5214837}{679601578} & \frac{424549711}{451760648} \\
 -\frac{153720698}{116964597} & \frac{98271451}{22705510} & \frac{108839271}{102671668} & -\frac{37652132}{76331505} & -\frac{98556781}{98719039} & -\frac{86545565}{156277133} \\
 -\frac{514694857}{175951034} & \frac{108839271}{102671668} & \frac{178543136}{16280101} & -\frac{13588975}{554452603} & -\frac{66382289}{197496474} & -\frac{315010733}{72953806} \\
 -\frac{26958123}{134065612} & -\frac{37652132}{76331505} & -\frac{13588975}{554452603} & \frac{127914572}{23010911} & -\frac{93779957}{771873704} & -\frac{258311971}{316622401} \\
 -\frac{5214837}{679601578} & -\frac{98556781}{98719039} & -\frac{66382289}{197496474} & -\frac{93779957}{771873704} & \frac{183401329}{33290110} & -\frac{60904303}{161208591} \\
 \frac{424549711}{451760648} & -\frac{86545565}{156277133} & -\frac{315010733}{72953806} & -\frac{258311971}{316622401} & -\frac{60904303}{161208591} & \frac{502508117}{78490640} \\
\end{array}
\right) \succeq 0\\
Q_2 &= \left(
\begin{array}{cccccc}
 \frac{21520940}{25577879} & -\frac{56020343}{32074003} & \frac{40731578}{751516279} & -\frac{46544963}{139367268} & -\frac{41177990}{108764983} & -\frac{26606007}{46612636} \\
 -\frac{56020343}{32074003} & \frac{112841678}{19842961} & -\frac{139240153}{172670104} & -\frac{45501317}{43903809} & \frac{64055491}{88725341} & \frac{21288583}{30121110} \\
 \frac{40731578}{751516279} & -\frac{139240153}{172670104} & \frac{68362401}{21097442} & -\frac{168386141}{819717774} & -\frac{155286027}{198655888} & \frac{30506956}{19158511} \\
 -\frac{46544963}{139367268} & -\frac{45501317}{43903809} & -\frac{168386141}{819717774} & \frac{166235485}{28138938} & \frac{15677552}{218059291} & -\frac{15992364}{25871383} \\
 -\frac{41177990}{108764983} & \frac{64055491}{88725341} & -\frac{155286027}{198655888} & \frac{15677552}{218059291} & \frac{253525900}{46511459} & \frac{95613053}{837681775} \\
 -\frac{26606007}{46612636} & \frac{21288583}{30121110} & \frac{30506956}{19158511} & -\frac{15992364}{25871383} & \frac{95613053}{837681775} & \frac{267687310}{41812157} \\
\end{array}
\right) \succeq 0\\
Q_3 &= \left(
\begin{array}{cccccc}
 \frac{29877454}{113194375} & -\frac{110018062}{390364861} & -\frac{39492021}{93889856} & -\frac{44736353}{260223501} & -\frac{27286543}{148218452} & -\frac{211317628}{549271497} \\
 -\frac{110018062}{390364861} & \frac{168343502}{34876437} & -\frac{63781869}{56201314} & \frac{813722845}{556876698} & \frac{1719950}{4084346189} & \frac{20149420}{711586093} \\
 -\frac{39492021}{93889856} & -\frac{63781869}{56201314} & \frac{293980380}{89098241} & -\frac{16659683}{50114131} & \frac{24295714}{57792167} & \frac{20062513}{28511329} \\
 -\frac{44736353}{260223501} & \frac{813722845}{556876698} & -\frac{16659683}{50114131} & \frac{166235485}{28138938} & -\frac{91733513}{894919007} & -\frac{11949058}{15299253} \\
 -\frac{27286543}{148218452} & \frac{1719950}{4084346189} & \frac{24295714}{57792167} & -\frac{91733513}{894919007} & \frac{187073509}{34708874} & -\frac{1990762}{36615949} \\
 -\frac{211317628}{549271497} & \frac{20149420}{711586093} & \frac{20062513}{28511329} & -\frac{11949058}{15299253} & -\frac{1990762}{36615949} & \frac{192697280}{35564393} \\
\end{array}
\right) \succeq 0.
\end{align*}

At this stage it does not matter (for the purposes of a verifiable proof) where $R_0,Q_1,Q_2,Q_3$ came from; it suffices to show that they are positive semidefinite and that the polynomial
$$
\eps(y,z) := -h(y,z) - zs_1 - (y-z)s_2 - (\alpha-y)s_3 - r_0
$$
is positive on $[0,1]^2$, where $r_0 = \underline{x}^\intercal R_0\underline{x}$. 
To check positive semidefiniteness of a matrix $A=(a_{ij})_{i,j\in[m]}$, we first check that $A$ is symmetric, then we use Sylvester's criterion, which says that a Hermitian matrix $A$ is positive semidefinite if and only if $A^{(k)}=(a_{ij})_{i,j\in[k]}$ has positive determinant for all $k \in [m]$.
We bound $\eps(y,z)$ from below by its constant term minus the sum of the absolute value of its other coefficients (see \texttt{311.nb}) to see that $\eps(y,z) \geq \frac{1}{50}$ in the required region.
This completes the proof of the claim.
\ecpf

Since Claim~\ref{cl:h} implies that $\OPT=\{(\frac{3}{5},0,\ldots)\}$,
it remains to check that $p(K_{3,1,1},\cdot)$ is strict.
Consider $G=G_{n,\V a}$ which has a clique part $V_0$ of size $\frac{2n}{5} + O(1)$ and another part $V_1$ of size $\frac{3n}{5}+O(1)$ which is an independent set. Now~\Sone~is immediate as $G\oplus xy$ has no induced copy of $F$ containing both $x$ and $y$.

Now we check~\Stwo. 
Let $c:=\frac{108}{125}$. 
We have $\supp^*(\V a)=\{0,1\}$, so given any $b: \{1\} \to \{0,1\}$ and $\alpha \in [0,1]$ it is enough to show that
$$
\join_{b,\alpha}\lambda(\V a) = \lambda(\V a) - \lambda(\V a,(b,\alpha)) = \lambda(\V a) -\lim_{n\to\infty}\lambda(G_{n,\V a} +_{b,\alpha}u,u) \geq \frac{2}{5}c(1-\alpha) = \lambdamax(1-\alpha),
$$ 
that is, $\lambda(\V a, (b,\alpha)) \leq \lambdamax \alpha$.
If $b(1)=0$ then $u$ lies in a copy of $K_{3,1,1}$ only if it lies in the $3$-set with two vertices in $V_1$ and the two singletons are in $N(u) \cap V_0$, so $\lambda(\V a,(b,\alpha)) = \binom{4}{2,2}\left(\frac{2\alpha}{5}\right)^2\left(\frac{3}{5}\right)^2 = \lambdamax \alpha^2$, as required.
If $b(1)=1$ then $u$ lies in a copy of $K_{3,1,1}$ only if the $3$-set is in $V_1$ and the other singleton is in $N(u) \cap V_0$, so
$\lambda(\V a,(b,\alpha)) = \binom{4}{1,3}\left(\frac{2\alpha}{5}\right)\left(\frac{3}{5}\right)^3 = \lambdamax \alpha$, as required.

This completes the proof of the theorem.
\epf

\section{Concluding remarks}\label{sec-conclude}

In this paper we have shown how to obtain stability from results in extremal graph theory which use symmetrisation. We have applied our general theory to the inducibility problem for complete partite graphs.
It would be interesting to solve other instances of the polynomial optimisation problem which amounts to determining $i(F)$. 

It would be particularly interesting to find other extremal graph theory problems to which our theory applies.

\bibliography{symm}
\end{document}